\documentclass[11pt,letterpaper,oneside,usenames]{article}
\usepackage{geometry}
\usepackage{caption}
\usepackage{subcaption}

\let\oldmarginpar\marginpar
\renewcommand\marginpar[1]{\-\oldmarginpar[\raggedleft\footnotesize #1]%
{\raggedright\footnotesize #1}}

\usepackage{amsmath}
\allowdisplaybreaks[1]
\usepackage{appendix}
\usepackage{lscape}
\usepackage{multirow}
\usepackage[stable]{footmisc}

\usepackage[a4,center,noinfo,dvips]{crop}

\usepackage[english]{babel}
\usepackage[OT1]{fontenc}
\usepackage{graphicx}
\usepackage{psfrag}
\usepackage{amssymb}
\usepackage{amsfonts}
\usepackage{tabularx}
\usepackage{fancybox}
\usepackage{setspace}
\usepackage{hhline}
\usepackage[dvipsnames]{xcolor}
\usepackage{makeidx}
\usepackage{enumerate}

\usepackage{longtable}
\usepackage{booktabs}

\makeindex

\usepackage{natbib}
\bibpunct{(}{)}{;}{a}{}{,}
\DeclareGraphicsExtensions{.ps,.eps}

\newtheorem{lemma}{Lemma}[section]
\newtheorem{prop}{Proposition}[section]

{
\begin{bmatrix}}%
{\end{bmatrix}
}

\newcommand{\myRed}[1]{{\color{BrickRed}{#1}}}

\newcommand{\ph}[1]{\phantom{#1}}
\newcommand{\uw}{\underline{w}}
\newcommand{\ow}{\overline{w}}
\newcommand{\mcA}{\mathcal{A}}
\newcommand{\mcN}{\mathcal{N}}
\newcommand{\mcL}{\mathcal{L}}
\newcommand{\mcH}{\mathcal{H}}

\usepackage{setspace}
\title{\textbf{The fragility-constrained vehicle routing problem with time windows}}

\usepackage{authblk} 

\newcommand{\corresp}{\thanks{Corresponding author \\
Email addresses:  \texttt{clement.altman@polytechnique.org} (Cl\'ement Altman), \\ 
\texttt{guy.desaulniers@gerad.ca} (Guy Desaulniers), \texttt{fausto.errico@etsmtl.ca} (Fausto Errico)}}

\author[1]{Cl\'ement Altman}
\author[1]{Guy Desaulniers}
\author[2]{Fausto Errico\corresp }

\affil[1]{{\small Department of Mathematics and Industrial Engineering, Polytechnique Montr{\'e}al, \& GERAD, Canada}}

\affil[2]{{\small Departement of Civil Engineering, \'{E}cole de Technologie Sup\'{e}rieure de Montr\'{e}al, GERAD \& CIRRELT, Canada}}

\begin{document}


\sloppy \thispagestyle{empty}

\newpage
\maketitle

\paragraph{Abstract:} 

We study a new variant of the well-studied Vehicle Routing Problem with Time Windows (VRPTW), 
called the fragility-constrained VRPTW, which 
assumes that 1) the capacity of a vehicle is organized in multiple identical stacks; 2) all items picked up at a customer 
are either ``fragile" or not; 3) no non-fragile items can be put on top of a fragile item 
(the \emph{fragility} constraint); and 4) no en-route load rearrangement is possible. 
We first characterize the feasibility of a route with respect to this fragility constraint. 
Then, to solve this new problem, we develop an exact branch-price-and-cut (BPC) algorithm that 
includes a labeling algorithm exploiting this feasibility characterization to 
efficiently generate feasible routes. This algorithm is benchmarked against 
another BPC algorithm that deals with the fragility constraint in the 
column generation master problem through infeasible path cuts. Our computational 
results show that the former BPC algorithm clearly outperforms the latter in 
terms of computational time and that the fragility constraint has a greater impact on 
the optimal solution cost (compared to that of the VRPTW) 
when vehicle capacity decreases, stack height increases and 
for a more balance mix of customers with fragile and non-fragile items. 

\vspace{-0.5cm}


\pagenumbering{arabic}
\section{Introduction}\label{sec:intro}

This paper introduces a new variant of the classical Vehicle Routing Problem with Time
Windows \citep[VRPTW, see][]{desaulniers+mr2014} where a particular \emph{fragility} 
loading constraint is enforced. This variant is called the Fragility-constrained VRPTW (F-VRPTW). The 
classical VRPTW considers a set of customers, each specifying a demand volume to be picked up 
and a pickup time window. Pickups are performed by a fleet of homogeneous vehicles with given capacity, initially located at a depot.
Travel times and costs among customers and between each customer and the depot are assumed to be known.
The VRPTW consists in computing least-cost vehicle routes such that each customer is visited exactly once within the indicated time window, while the vehicle capacity is not exceed for each route. It has been extensively studied in the literature \citep[see][]{desaulniers+mr2014} and the current state-of-the-art exact methodology for solving it is branch-price-and-cut \citep[BPC, see][]{costa+cd2019}. In particular, the  sophisticated BPC algorithms of \citet{pecin+cdu2017} and \citet{sadykov+up2020} are able to solve to proven optimality most tested instances with up to 200 customers. 

For the F-VRPTW, the vehicle capacity is organized in multiple identical stacks of a given maximal height. 
Customers might have items labeled as ``fragile''.
To avoid potential damages to these items, the F-VRPTW imposes the fragility
constraint, which implies the following:
1) non-fragile items are prevented to be stacked on top of fragile ones,
however
2) fragile items are allowed to be stacked on top of other
fragile items.
Furthermore, no en-route load rearrangement is allowed. 
This problem is cast as a pickup one but the developed 
methodology can easily be adapted for a delivery problem, 
simply by inverting the role of the fragile and non-fragile items. 
Our motivation in studying the F-VRPTW
is its relevance for naval shipping
companies delivering freight in northern Quebec and, to the best of our knowledge, this problem has never been studied before.
However, although
in a different loading framework (three-dimensional loading), the fragility constraint, in the exact sense used in the present work,
was already introduced by \cite{gendreau+ilm2006}.
Furthermore, the F-VRPTW is related to a number of studies
addressing routing problems combined with several types of loading constraints.
For an extensive literature review, the reader is referred to
\cite{pollaris+bcjl2015} and references therein.

A first stream of work related to the F-VRPTW focuses on routing problems where
the vehicle capacity is organized in multiple stacks \citep[see][to cite a few]{cote+asgp2012,carrabs+cs2013,cherkesly+dil2016,venstra+cdl2017}.
Handling multiple stacks responds to practical concerns in a variety of contexts,
such as grocery, pharmaceutical and naval applications
where either multi-compartments vehicles are required,
or freight need to be organized into pallets or containers
which are then stacked for shipping. From the operational viewpoint,
structuring the load in stacks limits the access to the only items
located on the top of the stacks, unless time-consuming rearrangements of the freight are operated.
In real-world operations, this limitation
results in the so-called LIFO policy.
Typically, the LIFO policy is a concern when the underlying routing decisions belong to the family of the Pickup and Delivery Problems (PDPs).
In fact, the pickup order will influence which items must be delivered first, thus deeply impacting routing possibilities.
In the VRPTW literature, however, the LIFO policy is generally not a concern
because all the freight either originates from the depot (in case of deliveries) where it can be well positioned 
in the vehicle before departure or is destined to the depot (in case of pickups) where it can be typically unloaded in any order. 
Due to the structural differences between the PDPs and the VRPTW,
providing an extensive literature review on the PDPs is out of the scope,
and we again refer the interested reader to \cite{pollaris+bcjl2015}. However, we notice that,
given the additional fragility constraint in the F-VRPTW,
the LIFO rule may render infeasible some customer sequences, and this needs to be taken into account when solving the F-VRPTW.

Another related stream of literature addresses the \emph{load-bearing strength} constraint,
which can be seen, to some extent, as a generalization
of the fragility constraint introduced here. This constraint assures that the total weight
stacked on top of a given item does not exceed a maximal
value, thus preventing item damages. \cite{junqueira+ocm2013}
proposes an integer linear programming model for the
vehicle routing problem with three-dimensional
loading constraints (3L-CVRP), including the load-bearing
strength. By adopting a commercial mixed-integer programming solver, the authors obtained the solution
for some small problem instances.
More recently, a single-stack
variant of this problem was addressed by \cite{chabot+lcr2016}
in the context of warehouse order picking activities.
The authors propose several heuristics as well as
exact branch-and-cut algorithms.
As previously mentioned, \cite{gendreau+ilm2006} address the 3L-CVRP where
the fragility constraint is handled as in our setting.
For this problem, the authors developed an efficient Tabu Search (TS) algorithm.
Different metaheuristics have been developed for the same
problem by several other authors, see
\cite{fuellerer+dhi2010,tarantilis+zk2009,tao+w2015}, to cite a few.

In this paper, we propose a set partitioning formulation of the F-VRPTW, which we solve
via a BPC algorithm. The most challenging aspect is in the
solution of the column generation pricing problem, which
turns out to be a variant of the Elementary Shortest Path Problem with Resource Constraints
(ESPPRC). The main contribution of this paper is
in the way we handle the fragility constraint.
We solve the pricing problem by a labeling algorithm and
provide a formal characterization of partial routes allowing us to
easily check their feasibility with respect to the fragility constraint.
In particular, this characterization
allows us to avoid the intuitive but naive approach of duplicating labels to account for alternative freight configurations. To improve the performance of the labeling algorithm,
we develop an efficient dominance rule with the purpose of discarding proven dominated labels.

To analyze the performance of our algorithm, we perform an extensive computational campaign. We first introduce F-VRPTW instances that are obtained from the Solomon's VRPTW benchmark instances by varying the vehicle capacity and the maximal stack height.
Our computational results show that the fragility constraint is more binding when vehicle capacity is smaller (i.e., for shorter routes) and when stacks are higher. 
Furthermore, to assess the performance of the proposed BPC algorithm,
we have developed an alternative solution method, which is a classic BPC algorithm for the VRPTW (possibly generating routes violating the fragility constraint), augmented with infeasible path cuts to eliminate integer solutions violating the
fragility constraint.
Our test results show that the proposed BPC algorithm consistently outperforms the alternative one.

The rest of the paper is organized as follows. In Section \ref{sec:formulation}, we formally state the F-VRPTW and
present a mathematical formulation. In Section \ref{sec:char}, we prove the main theoretical result of the paper,
that gives a formal characterization of a route feasible with respect to the fragility constraint.
This result is then exploited in Section \ref{sec:BPC} to develop an efficient BPC algorithm.
We report and analyze computational results in Section \ref{sec:experiments}, before providing concluding remarks in Section~\ref{sec:conclusions}.

\section{Problem statement and mathematical formulation}\label{sec:formulation}
The classical VRPTW can be described as follows. Consider a directed graph $G=(\mcN, \mcA)$ where
$\mcN=\{0,1,\ldots,n, n+1\}$ is the node set and $\mcA=\{(i,j)~|~i,j \in N,~i\neq j, (i,j) \neq (n+1, 0)\}$ is the arc set.
We denote $\mcN_c=\{1,\ldots,n\}$ the subset of nodes representing the $n$ customers, while nodes $0$ and 
$n+1$ represent the same depot at the beginning and the end of a route, respectively.
A fleet of homogeneous vehicles with capacity Q is initially located at the depot.
A demand volume $q_i$ (an integer number of items of the same size) 
to be collected and a time window $[\uw_i,\ow_i]$ are
associated with each customer $i\in \mcN_c$. A vehicle is allowed to arrive at customer $i$
earlier than $\uw_i$, but must leave in the specified time window.
A cost $c_{ij}$ and a travel time $t_{ij}$ (possibly including a service time at $i$) are associated with each arc $(i,j)\in \mcA$. We assume that these costs and travel times satisfy the triangle inequality. Observe 
that set $\mcA$ can be reduced by eliminating all arcs $(i,j)$ with $i,j\in \mcN_c$ that are 
load or time infeasible, i.e., such that $q_i+q_j > Q$ or $\uw_i + t_{ij} > \ow_j$. 
The VRPTW calls for finding a minimum-cost set of routes starting and ending at the depot,
such that each customer is visited exactly once, and the time windows and vehicle capacity constraints are fulfilled. 

\begin{figure}[t]
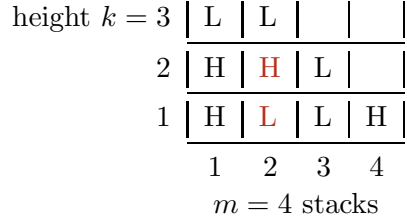

\begin{center}
\begin{tabular}{r|l|l|l|l|l|}
 height $k=3$ & L  & L &  & \\ \cmidrule{2-5}
 2 &  H & \myRed{H} & L &  \\ \cmidrule{2-5}
 1 & H & \myRed{L}  &  L & H  \\ \cmidrule{2-5}
 \multicolumn{1}{c}{} & \multicolumn{1}{c}{1} & \multicolumn{1}{c}{2} & \multicolumn{1}{c}{3} 
 & \multicolumn{1}{c}{4} \\
 \multicolumn{1}{c}{} & \multicolumn{4}{c}{$m=4$ stacks}
\end{tabular}
\end{center}
\caption{An infeasible loading configuration for a vehicle with 4 stacks of height 3}
\label{fig:loads}
\end{figure}

The F-VRPTW differs from the VRPTW in that the vehicle capacity is organized
in a set of $m$ identical stacks $M=\{1,\ldots,m\}$, each with a set of 
positions $K = \{ 1, \ldots, k\}$ (see Figure \ref{fig:loads}). Consequently, $k$ is the height of 
a stack and the vehicle 
capacity is set as $Q=km$. Items are assumed to be put in the vehicle from the top 
and each of them occupies one position in one stack.
Furthermore, the customer set $\mcN_c$ is partitioned in two
subsets $\mcL$ and $\mcH$, identifying customers with fragile (light) and non-fragile (heavy) items,
respectively.
The fragility constraint requires that
no item from customers in $\mcH$ can be stacked on the top of items from customers
in $\mcL$. However,
fragile items can be stacked on top of other fragile ones. For example,
the positioning of the items in stack 2 in Figure \ref{fig:loads} is infeasible, while
the positioning in the other stacks is legitimate. No
en-route rearrangement of items is permitted. It should be noticed
that the case of customers having both fragile and non-fragile items
can be easily reduced to the present setting by adding
suitable dummy nodes. \\

We formulate the F-VRPTW as a set-partitioning problem. To this scope,
consider a route $r$ as a sequence of nodes $r=(v_0, v_1,\ldots, v_p, v_{p+1})$, where
$v_0 = 0$ and $v_{p+1} = n+1$ represent the depot, and let $\Omega$ be set of all feasible routes.
For all $i\in \mcN_c$, let $a_{ir}$ be a parameter with value 1 if route $r$ visits customer $i$ and $0$ otherwise.
The cost $c_r$ of route $r\in \Omega$ is computed as
\begin{equation}\label{eq:routeCost}
c_r=\sum_{i=0}^p c_{v_i,v_{i+1}}.
\end{equation}
By introducing a binary variable $x_r$ for each route $r\in \Omega$ that takes value 1 if route $r$ is chosen and 0 otherwise, the F-VRPTW can be formulated as follows:
\begin{align}
\min \quad & \sum_{r\in \Omega} c_r x_r															 & \label{eq:costs}\\
s.t. \quad & \sum_{r \in \Omega} a_{ir}x_r = 1,										 & \forall i \in \mcN_c \label{eq:covering}\\
		 & x_r \in \{0,1\} ,																					 & \forall r \in \Omega \label{eq:binary},
\end{align}
where the objective function \eqref{eq:costs} minimizes the total cost, constraints \eqref{eq:covering} ensure
that each customer is visited exactly once, and constraints \eqref{eq:binary} impose that the variables are  
binary.

\section{Characterization of fragility-feasible routes}\label{sec:char}

The present section is devoted to
introducing a formal characterization of a route that is feasible with respect to
the fragility constraint (hereafter called a fragility-feasible route). As previously mentioned,
the proposed solution method heavily stands on this characterization.
In particular, as it will be clear
in Section \ref{sec:BPC}, this characterization
allows us to drastically reduce the number of labels generated
in the labeling algorithm developed to solve the column generation pricing problem.

The main observation behind the proposed characterization is that
a route can be infeasible with respect to the fragility constraint only if, in the
collecting process,
we must form an incomplete stack containing fragile items, and then we keep collecting
a large enough number of non-fragile items so that some of them must be placed on the incomplete
stack with fragile items, hence violating the fragility constraint. As a consequence,
the number of non-fragile items that can be collected after a given customer is bounded by a function depending
on the number of fragile and non-fragile items collected so far and the order in which they were collected. However, when 
it is possible to complete all the stacks containing fragile items at a given customer with fragile items, the number of non-fragile
items that can still be collected is only bounded by the vehicle capacity constraint.

Given a route $r=(v_0, v_1,\ldots, v_p, v_{p+1})$, let us consider the induced sequence of individual collected items
$S_r=(b_1,b_2, \ldots, b_s)$.
Given that customers belong either to set $\mcL$ or $\mcH$, the order used to embed the items
of a given customer in sequence $S_r$ is irrelevant. We denote by $a^\mcL(i)$  the number of
fragile items in $S_r$ until and including item $b_i$, by $a^\mcH(i)$ the corresponding number of 
non-fragile items, and by $l(i)=a^\mcL(i)+a^\mcH(i)$ the total number of collected items. 
Even if $l(i) = i$ here, notation $l(i)$ is used because the items will not always be considered individually 
in the rest of this paper. Let us also define a modified modulo function as:
\begin{equation} F(x) =
\left \{
\begin{array}{l}
    x \bmod k  \hspace{0.4cm} \text{if} \hspace{0.2cm} x \bmod k \hspace{0.2cm} \neq \hspace{0.2cm} 0, \\
    k   \hspace{1.5cm} \text{otherwise}. \\
\end{array}
\right.
\end{equation}
Finally, by a slight abuse of notation, we will write $b_i \in \mcL$ or $b_i \in \mcH$ to express
that $b_i$ is a fragile or non-fragile item, respectively.

\begin{prop}[\bf Characterization]
\label{prop:characterization}

Consider a route $r$ and let us assume that $r$ is feasible
with respect to the capacity constraint. 
The corresponding sequence $S_r=(b_1,b_2, \ldots, b_s)$ is feasible 
with respect to the fragility constraint
if and only if, for all $i\in \{ 1, \ldots, s \}$ with $b_i \in \mcL$, either
\begin{enumerate}[(c1)]
\item $a^\mcH(i) + F(a^\mcL(i)) \geq k$
\end{enumerate}
or
\begin{enumerate}[(c2)]
\item $a^\mcH(s) - a^\mcH(i) \leq U(i)$,
\end{enumerate}
where $a^\mcH(s)$ denotes the total number of non-fragile items and $U(i)= mk - (l(i)-l(i)\bmod k +k)$ is an upper bound on the number of non-fragile items that can be collected after item~$b_i$.
\end{prop}

Before proving this proposition, let us present the intuition behind conditions  $(c1)$ and $(c2)$, and how they can be interpreted. 
Each fragile item $b_i\in \mcL$ in $S_r$ may induce an upper bound on the number of 
non-fragile items that can be collected after it. This upper bound corresponds to the optimal value of 
the problem of finding a feasible loading configuration for the first $i$ items in $S_r$ 
that maximizes the number of positions available to load additional non-fragile 
items. Note that this upper bound does not depend on the items to be collected after item $b_i$ and the corresponding configuration is not necessarily a sub-configuration of a feasible configuration for the whole sequence $S_r$. Consequently, this configuration is called, hereafter, a locally-optimal configuration for item $b_i$ (even if $b_i \in \mcH$). 
When condition $(c1)$ is not met for $b_i$, i.e., $a^\mcH(i) + F(a^\mcL(i)) < k$, 
then there remain in any locally-optimal configuration for $b_i$ exactly $k - l(i) \bmod k$ 
empty positions above a fragile item to which no non-fragile item can be assigned. 
Therefore, the maximum number of non-fragile items that can be loaded after item $b_i$ is equal to 
$mk - l(i) - (k - l(i) \bmod k) = U(i)$, where $mk - l(i)$ is the total number of empty positions. 
On the other hand, when condition $(c1)$ holds for $b_i$, we can show 
that either 1) item $b_i$ has no impact on the number of additional non-fragile items that 
can be loaded, or 2) the corresponding upper bound is equal to $mk - l(i)$ and, thus, 
redundant with the vehicle capacity constraint. In both cases, there is no need 
to compute the upper bound. 

To illustrate these different cases, let us consider the following example. Let 
$m = 4$, $k = 3$, and $S_r = (b^\mcL_1,b^\mcH_2, b^\mcL_3, b^\mcL_4, b^\mcL_5, b^\mcL_6, 
b^\mcH_7, b^\mcL_8, b^\mcL_9)$, where the upper index indicated if the item belongs to $\mcL$ or $\mcH$. 
Table \ref{tab:example} displays relevant values for each item $b^\mcL_i$: 
the left-hand side of condition $(c1)$, i.e., $a^\mcH(i) + F(a^\mcL(i))$; 
the upper bound $U(i)$ of condition $(c2)$ whenever condition 
$(c1)$ does not hold; the upper bound $mk-l(i)$ induced by the 
vehicle capacity; and the exact maximum number of non-fragile items 
that can be loaded after item $b^\mcL_i$ ($max^\mcH(i)$). The latter is computed as: 
\begin{equation}
max^\mcH(i) = \min \> \{ mk - l(i), \min\limits_{j \in B_i } \{ U(j) - (a^\mcH(i) - a^\mcH(j)) \}  \}, 
\end{equation}
where $B_i$ is the subset of fragile item indices $j \in \{ 1, \ldots, i\}$ such that $a^\mcH(j) + F(a^\mcL(j)) \ge k$. 
This maximum is obtained by considering the upper bounds induced by the vehicle capacity 
and by the fragile items $b^\mcL_j$, $j < i$, for which condition $(c1)$ does not hold. In the latter case, 
the corresponding upper bound $U(j)$ is adjusted by subtracting the number of non-fragile items 
loaded between $b^\mcL_j$ and $b^\mcL_i$. 

\begin{table}[t]
\centering 
\scalebox{0.8}{
\begin{tabular}{c|cccc}
\toprule
$i$ & $a^\mcH(i) + F(a^\mcL(i))$ & $U(i)$ & $mk - l(i)$ & $max^\mcH(i)$ \\
\midrule
1 & 1 & 9 & 11 & 9 \\
3 & 3 & - & 9 & 8 \\
4 & 4 & - & 8 & 8 \\
5 & 2 & 6 & 7 & 6 \\
6 & 3 & - & 6 & 6 \\
8 & 5 & - & 4 & 4 \\
9 & 3 & - & 3 & 3 \\
\bottomrule
\end{tabular}
}
\caption{Values for items in $\mcL$ for an example with $m=4$ and $k=3$}
\label{tab:example}
\end{table}

\begin{figure}[t]
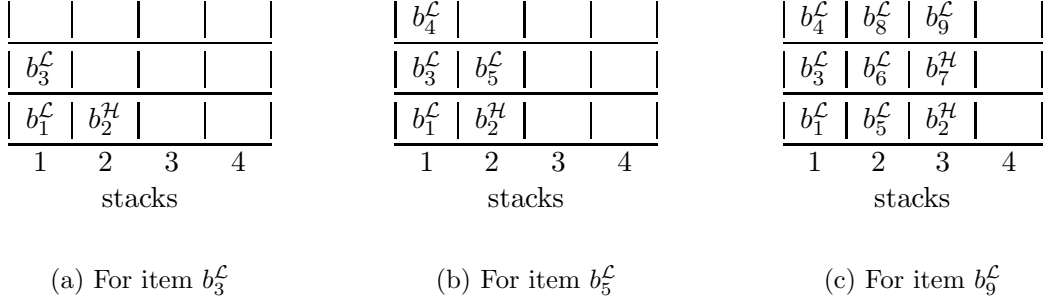

\centering 
\begin{subfigure}[t]{.3\textwidth}
\begin{center}
\begin{tabular}{|l|l|l|l|l|}
   &  &  & \\ \midrule
 $b^\mcL_3$ &  &  &  \\ \midrule
 $b^\mcL_1$ & $b^\mcH_2$ & \ph{$b^\mcH_7$} & \ph{$b^\mcH_7$}  \\ \bottomrule
\multicolumn{1}{c}{1} & \multicolumn{1}{c}{2} & \multicolumn{1}{c}{3} 
 & \multicolumn{1}{c}{4} \\
 \multicolumn{4}{c}{stacks}
\end{tabular}
\end{center}
\captionof{figure}{For item $b^\mcL_3$}\label{fig:exL3}
\end{subfigure}
\begin{subfigure}[t]{.3\textwidth}
\begin{center}
\begin{tabular}{|l|l|l|l|l|}
 $b^\mcL_4$  &  &  & \\ \midrule
 $b^\mcL_3$ &  $b^\mcL_5$ &  &  \\ \midrule
 $b^\mcL_1$ & $b^\mcH_2$ & \ph{$b^\mcH_7$} & \ph{$b^\mcH_7$}  \\ \bottomrule
\multicolumn{1}{c}{1} & \multicolumn{1}{c}{2} & \multicolumn{1}{c}{3} 
 & \multicolumn{1}{c}{4} \\
 \multicolumn{4}{c}{stacks}
\end{tabular}
\end{center}
\captionof{figure}{For item $b^\mcL_5$}\label{fig:exL5}
\end{subfigure}
\begin{subfigure}[t]{.3\textwidth}
\begin{center}
\begin{tabular}{|l|l|l|l|l|}
 $b^\mcL_4$  & $b^\mcL_8$ & $b^\mcL_9$ & \\ \midrule
 $b^\mcL_3$ &  $b^\mcL_6$ & $b^\mcH_7$ &  \\ \midrule
 $b^\mcL_1$ & $b^\mcL_5$ & $b^\mcH_2$ & \ph{$b^\mcH_7$}  \\ \bottomrule
\multicolumn{1}{c}{1} & \multicolumn{1}{c}{2} & \multicolumn{1}{c}{3} 
 & \multicolumn{1}{c}{4} \\
 \multicolumn{4}{c}{stacks}
\end{tabular}
\end{center}
\captionof{figure}{For item $b^\mcL_9$}\label{fig:exL9}
\end{subfigure}
\caption{Locally-optimal configurations for some items of the example}
\label{fig:example}
\end{figure}

Let us discuss the information reported in Table \ref{tab:example} for some specific items $b^\mcL_i$. 
For $i = 1$, the values are obvious because item $b^\mcL_1$ must be assigned 
to a stack, leaving the other three stacks (i.e., 9 positions) to potentially load non-fragile items. For $i=3$, 
condition $(c1)$ is satisfied ($a^\mcH(3) + F(a^\mcL(3)) = 3 \ge 3$). In this case, $b^\mcL_3$ can be 
positioned on top of a fragile item ($b^\mcL_1$) in the locally-optimal configuration for $b^\mcH_2$ 
(see Figure \ref{fig:exL3}) and, therefore, does not impact the number of non-fragile items 
that can still be loaded.  
Note that $max^\mcH(3) = 8$ stems from the upper bound $U(1) = 9$ from which item $b^\mcH_2$ 
has been subtracted. For $i = 5$, condition $(c1)$ is not satisfied ($a^\mcH(5) + F(a^\mcL(5)) = 2 < 3$) 
and the locally-optimal configuration (see Figure \ref{fig:exL5}) yields $U(5) = 6$ as at least one 
($= k - l(5) \bmod k$) of the seven ($= mk - l(5)$) unoccupied positions cannot be 
assigned to a non-fragile item (e.g., the empty position in stack 2 in Figure \ref{fig:exL5}). 
Note that item $b^\mcL_5$ induces a drop of two (from $max^\mcH(4) = 8$ to $max^\mcH(5) = 6$) 
in the exact maximum number of non-fragile items that can still be loaded. Finally, for $i = 9$, 
condition $(c1)$ holds. In this case, the locally-optimal configuration for $b^\mcL_9$ (see Figure \ref{fig:exL9}) 
provides an upper bound that is equal to $mk - l(9) = 3$ and is, thus, redundant with the 
vehicle capacity bound. 

To prove Proposition \ref{prop:characterization}, we proceed 
by induction on the number of stacks $m$. To make explicit the dependence
of the proposition on the number of stacks, we will use the notation $U(i,m)$ instead of $U(i)$.
The main line of the proof consists first in showing that the basis of the induction $P(m)$ with
 $m=1$ is true. In particular we show that the feasibility and characterization
conditions are equivalent: $Feas(1)\iff Char(1)$. We then prove the induction step $P(m) \implies P(m+1)$. In particular,
by assuming that the characterization is true for a generic number $m$ of stacks, i.e., $Feas(m)\iff Char(m)$, we show that $Feas(m+1) \impliedby Char(m+1)$
as well as $Feas(m+1)\implies Char(m+1)$, thus proving that $Feas(m+1)\iff Char(m+1)$ and completing the argument.

\subsection*{P(1) is true}

\noindent \emph{Proof.}

To prove that $P(1)$ is true, we need to show that $Feas(1)\iff Char(1)$.
Observe that the case with no fragile item is trivial because all routes are feasible with respect to the fragility constraint and,
in this case, the characterization requires no condition to be verified.
We can then concentrate on routes including fragile items ($a^\mcL(s) > 0$). We prove the sufficiency part first and the necessity second.
\begin{itemize}
\item $\bf Feas(1)\implies Char(1)$
Feasibility implies that non-fragile items are not stacked on top of fragile ones, i.e., 
$a^\mcH(s)-a^\mcH(i)=0$ for all $b_i \in \mcL$. For all $b_i \in \mcL$  such that $l(i)<k$, we have that $U(i,1)= 0$,
 and condition $(c2)$ is satisfied. Furthermore, if there exists $b_i \in \mcL$ such that $l(i) = k$, then 
 $F(a^\mcL(i)) = a^\mcL(i)$ and $a^\mcH(i) + F(a^\mcL(i)) = a^\mcH(i) + a^\mcL(i) = l(i) \ge k$ and condition (c1) is 
 verified. 
\item $\bf Feas(1)\impliedby Char(1)$
We prove the following counterpositive version of this statement:

\noindent \emph{If sequence $S_r$ is infeasible with respect to the fragility constraint, then
there exists an item $b_i \in \mcL$ such that
\begin{enumerate}[(c3)]
\item $a^\mcH(i) + F(a^\mcL(i)) < k$
\end{enumerate}
and
\begin{enumerate}[(c4)]
\setcounter{enumi}{1}
\item $a^\mcH(s) - a^\mcH(i) > U(i,1)$.
\end{enumerate}
}
If $S_r$ is infeasible, then a non-fragile item is collected after a fragile one and there exists $i \in \{ 1, \ldots, s \}$ such that $b_i \in \mcL$ and $a^\mcH(s) - a^\mcH(i) > 0$. Furthermore, $l(i) < k$ because an item is loaded after $b_i$. Consequently, $a^\mcH(i) + F(a^\mcL(i)) = a^\mcH(i) + a^\mcL(i) = l(i) <k$, thus condition $(c3)$ is satisfied. Also, $U(i,1) = 0$, and condition $(c4)$ follows. $\blacksquare$
\end{itemize}

\subsection*{P(m) $\implies$ P(m+1)}

\noindent \emph{Proof.}

Here we prove that $P(m+1)$ is true when assuming that $P(m)$ is true,
and in particular we need to prove that $Feas(m+1)\iff Char(m+1)$. We do this
by proving the necessity condition first and sufficiency second.

\noindent $\bf Feas(m+1)\impliedby Char(m+1)$

We start by considering a route
 $r$ collecting at most $(m+1)k$ items and complying with $Char(m+1)$. We build a suitable
subsequence ${\bar S}_r$ of at most $mk$ items and we show that ${\bar S}_r$ complies with 
$Char(m)$. By the induction hypothesis, this implies that ${\bar S}_r$ is a feasible
sequence. The feasibility of  $ S_r$ will be trivially obtained by construction.

We distinguish two cases:
\begin{enumerate}
\item $a^\mcL(s) \geq k$. In this case, route $r$ collects enough fragile items to be able to complete a stack
made of fragile items only. Consider the sequence ${\bar S}_r$ of at most $mk$ items obtained from $S_r$ by
eliminating the first
$k$ fragile items. Formally,
${\bar S}_r=(b_i~|~ b_i \in S_r \text{ and either } b_i \in \mcH \text{ or } b_i \in \mcL \text{ and } a^\mcL(i)>k)$.
We use bars to indicate that a given
entity refers to ${\bar S}_r$; for example, ${\bar a}_\mcH(i)$ is the number of non-fragile items collected
in ${\bar S}_r$ until and including item $b_i$.

Given that the number of non-fragile items collected in $S_r$ and ${\bar S}_r$ is the same and that 
exactly $k$ fragile items have been deleted from $S_r$, 
for all $i\in \{ 1, \ldots, s\}$ with $b_i \in {\bar S}_r$, we have 
\begin{eqnarray}
\label{eq:proof10}
{\bar a}^\mcH(i) & = &a^\mcH(i) \\
\label{eq:proof20}
{\bar a}^\mcL(i) & = & a^\mcL(i)-k \\
\label{eq:proof30}
{\bar l}(i) & = & l(i)-k.
\end{eqnarray}

To prove that ${\bar S}_r$ is feasible, we show that it satisfies $Char(m)$.
To this scope, let us consider a generic fragile item $b_i$ in ${\bar S}_r$ . 
Because we assume that $Char(m+1)$ holds for $S_r$, $b_i$ satisfies 
condition $(c1)$ or condition $(c2)$ for $S_r$, yielding the following two cases.
\begin{enumerate}[i)]
\item Condition $(c1)$ holds, i.e., $a^\mcH(i)+F((a^\mcL(i)))\geq k$. Equations \eqref{eq:proof10} and \eqref{eq:proof20}
imply $\bar{a}^\mcH(i)+F((\bar{a}^\mcL(i)))\geq k$, hence $b_i$ satisfies condition $(c1)$ for $\bar{S}_r$.
\item Condition $(c1)$ holds, i.e., $a^\mcH(s)-a^\mcH(i) \leq U(i, m+1)$. To prove the statement, we show that this implies  $\bar a^\mcH(s)-\bar a^\mcH(i) \leq \bar U(i, m)$, hence $b_i$ satisfies condition $(c2)$ for $\bar{S}_r$.
\begin{align}
\bar a^\mcH(s)-\bar a^\mcH(i) &= a^\mcH(s)-a^\mcH(i)  \label{eq:proof40} \\
& \leq  ~ U(i,m+1)   \label{eq:proof50} \\
&= (m+1)k - (l(i)-l(i)\bmod k +k)  \label{eq:proof60} \\
& =mk - (\bar l(i) - \bar l(i) \bmod k +k) \label{eq:proof70} \\
 &= \bar U(i,m),
\end{align}
where equations \eqref{eq:proof40} and \eqref{eq:proof70} hold 
given \eqref{eq:proof10} and \eqref{eq:proof30}, respectively.
\end{enumerate}
As a consequence of the induction hypothesis, the sequence $\bar{S}_r$ is feasible 
when considering $m$ stacks. Then, for $m+1$ stacks, $S_r$ is also 
feasible because it is always possible to stack the first $k$
fragile items (that were put aside when building $\bar{S}_r$) 
in a dedicated stack.
\item $a^\mcL(s) < k$.
We use a similar technique as for the previous case. Because it is not possible to fill
a stack with fragile items only, the strategy consists in building a dedicated stack
with a number of non-fragile items collected before the first fragile item, and with all 
fragile items. In case the number of non-fragile items collected before
the first fragile item plus the total number of fragile items is less than $k$,
the dedicated stack will be incomplete.

Let $f$ be the index in $S_r$ of the first fragile item. Consider the subsequence
$\bar S_r$ obtained from $S_r$ by eliminating the first
$min\{a^\mcH(f),k-a^\mcL(s)\}$ non-fragile items and $a^\mcL(s)$ fragile item. Because
$\bar S_r$ is only made of non-fragile items, it is certainly feasible with respect
to the fragility constraint. However, we need to ensure the capacity constraint
of $\bar S$, i.e., that $\bar l(s)\leq mk$ is fulfilled. We distinguish two cases:
\begin{itemize}
\item $k - a^\mcL(s)\leq a^\mcH(f)$. In this case, the construction eliminates exactly $k$ items and $\bar l(s)=l(s)-k \leq mk$.
\item $k - a^\mcL(s)> a^\mcH(f)$. In this case, the construction eliminates $a^\mcH(f)+a^\mcL(s)$ items. Hence,
\begin{align}
\bar l(s) & = l(s)-a^\mcH(f)-a^\mcL(s) \label{eq:proof80} \\
				& = a^\mcH(s)-a^\mcH(f) \label{eq:proof90} \\
        & \leq U(f,m+1) \label{eq:proof100} \\
				& = (m+1)k - (l(f)-l(f) \bmod k +k) \label{eq:proof110} \\
        & = mk \label{eq:proof120}.
\end{align}
Equation \eqref{eq:proof90} holds because $l(s)=a^\mcH(s)+a^\mcL(s)$. Inequality
\eqref{eq:proof100} is due to the hypothesis that $S_r$ complies with $Char(m+1)$ and, 
in particular, with condition $(c2)$ for $b_f$ because $k - a^\mcL(s) > a^\mcH(s)$ implies 
$a^\mcH(f) + a^\mcL(f) <k$ and, thus, that condition $(c1)$ ($a^\mcH(f) + F(a^\mcL(f)) \geq k$) does not hold. Finally, 
equation \eqref{eq:proof120} holds because $l(f) =  a^\mcH(f) + a^\mcL(f) < k$.
\end{itemize}
Consequently, $\bar S_r$ is feasible and $S_r$ is also feasible, because it is always possible
to build a dedicated stack with the eliminated items. This completes the proof of the necessity part of the statement.
$\blacksquare$
\end{enumerate}

\noindent $\bf Feas(m+1)\implies Char(m+1)$

We distinguish two cases:
\begin{enumerate}
\item $a^\mcL(s) \geq k$. 
Here, we use some of the constructions developed in the above necessity proof.
We start with a sequence $S_r$ of at most $(m+1)k$ items, which we assume to be feasible.
Then, we construct a suitable sequence $\bar S_r$ of at most $mk$
items and show that it is feasible. By using the induction hypothesis ($Feas(m) \iff Char(m)$),
we then prove that $S_r$ complies with $Char(m+1)$. 

As in the necessity proof, let $\bar S_r$ be the subsequence obtained by
eliminating the first $k$ fragile items from $S_r$.
Because we assume that $S_r$ is feasible and it is always
possible to put the first $k$ fragile items in a dedicated stack without affecting the feasibility of $S_r$,
the subsequence $\bar S_r$ is also
feasible and, by the induction hypothesis, $\bar S_r$ complies with $Char(m)$.
We also observe that relations \eqref{eq:proof10}-\eqref{eq:proof30} hold by construction.

To prove that $S_r$ complies with $Char(m+1)$, we first consider the fragile items $b_i$
belonging to both $S_r$ and $\bar S_r$ and observe that, by the same reasoning done for the necessity proof,
they all satisfy at least one of the conditions of $Char(m+1)$. In fact, for a given 
$i \in \{ 1, \ldots, s\}$ such that $b_i \in \mcL \cap \bar S_r$, we have two possible cases:
\begin{enumerate}[i)]
\item $\bar a^\mcH(i)+F(\bar a^\mcL(i))\geq k$. Then equations \eqref{eq:proof10} and \eqref{eq:proof20}
imply $a^\mcH(i)+F(a^\mcL(i)) \geq k$, hence $b_i$ satisfies condition $(c1)$ for $S_r$.
\item $\bar a^\mcH(s)-\bar a^\mcH(i) \leq \bar U(i,m)$. In this case, we get
\begin{align*}
a^\mcH(s)- a^\mcH(i) &= \bar a^\mcH(s)- \bar a^\mcH(i) \\
& \leq  ~  \bar U(i,m) \\
& =mk - (\bar l(i) - \bar l(i) \bmod k +k) \\
&= (m+1)k - (l(i)-l(i)\bmod k +k) \\
 &=  U(i,m+1).
\end{align*}
Hence, condition $(c2)$ holds for $S_r$.
\end{enumerate}

Now, let us verify that one of the conditions of $Char(m+1)$ also holds for the first $k$ fragile items belonging to $S_r$ only.
If $b_i$ is such that $a^\mcH(i) + F(a^\mcL(i)) \geq k$, condition~$(c1)$ of $Char(m+1)$ is
fulfilled. Otherwise, we must have $a^\mcH(i) + a^\mcL(i) < k$ (because $a^\mcL(i) \leq k$) 
and, consequently, $l(i) < k$. Thus,
\begin{equation}
a^\mcH(s)-a^\mcH(i) \leq mk = U(i,m+1),
\end{equation} where the inequality arises from the fact that $a^\mcL(s)\geq k$ and 
the equality holds because $l(i) < k$.

\item $a^\mcL(s) < k$. In this case we can prove the statement directly, without employing
the induction hypothesis. We have that $F(a^\mcL(i)) = a^\mcL(i)$ for all $i \in \{ 1, \ldots, s\}$ 
with $b_i\in \mcL$. Consequently, if item $b_i \in \mcL$ satisfies $a^\mcH(i) + a^\mcL(i) \geq k$, condition $(c1)$ of
$Char(m+1)$ is fulfilled. Otherwise, $a^\mcH(i) +a^\mcL(i) <k$, implying 
that $l(i) < k$ and $U(i,m+1)=mk$. In this case, condition $(c2)$ 
writes as $a^\mcH(s)-a^\mcH(i) \leq mk$. To prove that this condition is met, 
we start by observing that feasibility implies $a^\mcH(s) \leq mk + \min \{a^\mcH(f),k-a^\mcL(s)\}$, where $f$ is 
again the index of the first fragile item in $S_r$. 

We consider two cases:
\begin{enumerate}[i)]
\item $a^\mcH(f) < k-a^\mcL(s)$. In this case, $a^\mcH(s) \leq mk + a^\mcH(f)$ and we deduce 
\begin{align*} 
a^\mcH(s)-a^\mcH(i) & \leq a^\mcH(s)-a^\mcH(f) \\
& \leq  mk+a^\mcH(f)-a^\mcH(f) = mk.
\end{align*}
Therefore, condition $(c2)$ holds. 
\item $a^\mcH(f) \geq k-a^\mcL(s)$. In this case, $a^\mcH(s) \leq mk + k-a^\mcL(s)$ and we find
\begin{align*}
 a^\mcH(s)-a^\mcH(i) & \leq a^\mcH(s)-a^\mcH(f) \\
& \leq mk+k -a^\mcL(s) - a^\mcH(f) \\
& \leq mk,
\end{align*}
where the last inequality comes from the assumption that $a^\mcH(f) \geq k-a^\mcL(s)$, 
i.e., $k -a^\mcL(s) - a^\mcH(f) \le 0$. Consequently, condition $(c2)$ is also satisfied.
\end{enumerate}
\end{enumerate}
This completes the proof of the sufficiency part of the statement. $\blacksquare$

The characterization provided in Proposition \ref{prop:characterization} requires 
that conditions $(c1)$ and $(c2)$ are satisfied
for all fragile items $b_i$ in $S_r$. The following proposition ensures that,
if the characterization conditions are satisfied for a suitable subset of fragile items, they are satisfied for all
fragile items. In fact, it turns out that it is sufficient to check the conditions of the characterization only
for $b_i \in \mcL$ such that $i = s$ or $b_{i+1} \not\in \mcL$. In particular, when route $r$ visits a customer
$i \in \mcL$ with demand $q_i$, there is no need to consider these items individually. Moreover,
if $r$ visits several customers in $\mcL$ consecutively, it is sufficient that the last item of the last
of these customers verifies the conditions.

\begin{prop} \label{prop:simplerChar}
Let $i \in \{ 1, \ldots, s-1 \}$ and $b_i, b_{i+1} \in \mcL$. If 
\begin{equation*}
a^\mcH(i) + F(a^\mcL(i)) < k \quad \text{and} \quad a^\mcH(s) - a^\mcH(i) > U(i)
\end{equation*} 
(i.e., conditions $(c1)$ and $(c2)$ are not satisfied for $b_i$), then 
\begin{equation*}
a^\mcH(i+1) + F(a^\mcL(i+1)) < k \quad \text{and} \quad a^\mcH(s) - a^\mcH(i+1) > U(i+1)
\end{equation*} 
(i.e., they are neither satisfied for $b_{i+1}$). 
\end{prop}

\noindent\emph{Proof}

Let us denote the last two conditions $(c5)$ and $(c6)$. Observe that $a^\mcH(i+1) = a^\mcH(i)$, $a^\mcL(i+1) = a^\mcL(i) + 1$ and $l(i+1) = l(i)+1$. Furthermore, it is easy to prove that $U(i+1) \leq U(i)$. 

We get that $a^\mcH(s) - a^\mcH(i+1) = a^\mcH(s) - a^\mcH(i) > U(i) \geq U(i+1)$, which shows that condition $(c6)$ is always satisfied for $b_{i+1}$. 

Let us now consider two cases:
\begin{enumerate}
\item  If $a^\mcH(i) + F(a^\mcL(i)) < k - 1$, then 
\begin{align*}
a^\mcH(i+1) + F(a^\mcL(i+1)) & = a^\mcH(i) + F(a^\mcL(i) + 1) \\
& =  a^\mcH(i) + F(a^\mcL(i)) + 1 \\
& < k,
\end{align*}
where the second equality is valid because $F(a^\mcL(i)) < k - 1$. 
Thus, condition $(c5)$ is also satisfied.

\item If $a^\mcH(i) + F(a^\mcL(i)) = k - 1$, then $l(i) = a^\mcH(i) + a^\mcL(i) = k - 1 - F(a^\mcL(i)) + a^\mcL(i)$.
Because $a^\mcL(i) > 0$, we have that $a^\mcL(i) - F(a^\mcL(i)) = \psi k$, where $\psi$ is a nonnegative integer.
Thus, $l(i) = k - 1 + \psi k$, $l(i) \bmod k = k - 1$ and $U(i) = mk - l(i) - 1$.
The total number of items in the sequence is at least
$l(i+1) + a^\mcH(s) - a^\mcH(i+1) = l(i) + 1 + a^\mcH(s) - a^\mcH(i) > l(i) + 1 + U(i) = mk$,
which is not possible given the vehicle capacity.
Consequently, case 2 is not possible and the proof is complete. $\blacksquare$
\end{enumerate}

The labeling algorithm described in Section~\ref{sec:labeling}  
largely stands on Propositions \ref{prop:characterization} and~\ref{prop:simplerChar}. 

To conclude this section, we prove two other propositions concerning specific cases. 
Even if the stated results are not exploited in the proposed BPC algorithm, they are useful 
to understand when it is more difficult to satisfy the fragility constraint. Furthermore, they served 
to design our test instances and to analyze the obtained computational results. 

\begin{prop}\label{prop:stacks}
Consider a route $r = (v_0, v_1, \ldots, v_p, v_{p+1})$. 
If $\sum_{j=1}^p q_{v_j} \leq Q - k$, then $r$ satisfies the fragility constraint. 
\end{prop}

\noindent\emph{Proof}
Assume that $\sum_{j=1}^p q_{v_j} \leq Q - k$. From 
Propositions \ref{prop:characterization} and \ref{prop:simplerChar}, we can 
prove that $r$ is fragility-feasible by showing that, for the last item $b_{i_j}$ of 
every customer $v_j \in \mcL$, condition $(c2)$ is satisfied.  
Let $b_s$ be the last item of a customer $v_p$. If $j =p$, then $i_j = s$, $a^\mcH(s)- a^\mcH(i_j) = 0$, and $U(i_j,m) \ge 0$, implying condition $(c2)$. For $j < p$, we get
\begin{align*}
a^\mcH(s)- a^\mcH(i_j) & \leq \sum_{j'=j+1}^p q_{v_{j'}} \\
& \leq  Q - k - \sum_{j'=1}^{j} q_{v_{j'}} \\
& = mk - k - l(i_j) \\
& \leq mk - k - l(i_j) + l(i_j)\bmod k \\
 & =  U(i_j,m),
\end{align*}
where the first inequality comes from the fact that some customers $v_{j'}$, $j' \in \{ j+1, \ldots, p\}$, 
might belong to set $\mcL$ and the second inequality from the assumption. Therefore, 
condition $(c2)$ holds. 
$\blacksquare$

Proposition~\ref{prop:stacks} indicates that routes must be sufficiently filled to be fragility-infeasible. For example, 
if we consider $m=12$ stacks of height $k=4$, then only routes with a total load of 45 or more (i.e., with a filling rate exceeding 91.7\%) can violate the fragility constraint. Next, we discuss the case $k=2$.  

\begin{prop}\label{prop:k2}
Let $r = (v_0, v_1, \ldots, v_p, v_{p+1})$ 
be a route that is feasible with respect to the capacity constraint and let $k=2$. 
Route $r$ does not satisfy the fragility constraint if and only if $p \geq 2$, 
$\sum_{j=1}^p q_{v_j} = Q$ and there exists an index $j^* \in \{1, \ldots, p-1\}$ such that 
$v_j \in \mcL$ for all $j \in \{ 1, \ldots, j^*\}$, $v_j \in \mcH$ for all $j \in \{ j^*+1, \ldots, p\}$ and 
$\sum_{j=1}^{j^*} q_{v_j}$ is odd. 
\end{prop}

\noindent\emph{Proof}
($\implies$) Assume that $r$ does not satisfy the fragility constraint. 
From Propositions \ref{prop:characterization} and \ref{prop:simplerChar}, it means that there exists an 
index $j^*\in \{1, \ldots, p\}$ such that $v_{j^*} \in \mcL$ and both conditions $(c1)$ and $(c2)$ are not met for the 
last item $b_{i_{j^*}}$ of $v_{j^*}$. Because $k = 2$ and $F(a^\mcL(i_{j^*})) \in \{1, 2\}$, 
condition $(c1)$ is not satisfied if and only if $a^\mcH(i_{j^*}) = 0$ and $F(a^\mcL(i_{j^*})) = 1$. This means that 
$v_j \in \mcL$ for all $j \in \{ 1, \ldots, j^*\}$ and $\sum_{j=1}^{j^*} q_{v_j}$ is odd. Moreover, 
if $b_s$ denotes the last item of $v_p$, condition $(c2)$ writes as $a^\mcH(s) - 0 \le Q - (a^\mcL(i_{j^*}) - 1 + 2)$ 
because $l(i_{j^*}) = a^\mcL(i_{j^*})$ which is odd. This condition is thus violated if $a^\mcH(s) \geq Q - a^\mcL(i_{j^*})$, i.e., 
if $a^\mcH(s) + a^\mcL(i_{j^*}) = Q$ or, equivalently, $\sum_{j=1}^p q_{v_j} = Q$ and $v_j \in \mcH$ for all $j \in \{ j^*+1, \ldots, p\}$. 

($\impliedby$) Assume that $p \geq 2$, 
$\sum_{j=1}^p q_{v_j} = Q$ and there exists an index $j^* \in \{1, \ldots, p-1\}$ such that 
$v_j \in \mcL$ for all $j \in \{ 1, \ldots, j^*\}$, $v_j \in \mcH$ for all $j \in \{ j^*+1, \ldots, p\}$ and 
$\sum_{j=1}^{j^*} q_{v_j}$ is odd. Denote by $b_{i_{j^*}}$ the last item at customer $v_{j^*}$. Let us show that 
both conditions $(c1)$ and $(c2)$ of Proposition \ref{prop:characterization} do not hold for $b_{i_{j^*}}$, 
indicating that route $r$ is fragility-infeasible. From the assumption, we can easily deduce that 
$a^\mcH(i_{j^*}) = 0$ and $F(a^\mcL(i_{j^*})) = 1$. Thus, $(c1)$ is violated. Furthermore, we find that $a^\mcH(s) + a^\mcL(i_{j^*}) 
= \sum_{j=1}^p q_{v_j} = Q$ and $l(i_{j^*}) = a^\mcL(i_{j^*})$. Consequently, $a^\mcH(s) - a^\mcH(i_{j^*}) = Q - a^\mcL(i_{j^*}) > U(i_{j^*}) 
= Q - (a^\mcL(i_{j^*}) - 1 + 2)$, proving that condition $(c2)$ does not hold either. 
$\blacksquare$

Proposition \ref{prop:k2} describes explicitly the routes that do not satisfy the fragility constraint when $k=2$. 
In particular, it shows that, if a customer in $\mcH$ is visited before a customer in $\mcL$, then the route is necessary 
fragility-feasible. Indeed, it gives the opportunity to start a stack with a non-fragile item and fill it with a fragile item if needed. This principle, which is expressed by condition $(c1)$, is also valid for larger stack heights $k$. From it, we can deduce that routes visiting customers in $\mcH$ with sufficiently large demands before customers in $\mcL$ have high chances to be fragility-feasible. This may not be true for a large stack height $k$ but, in practice, $k$ is not too large. 


\section{The branch-price-and-cut framework}\label{sec:BPC}

BPC is a variant of the branch-and-bound algorithm
where the lower bounds at the nodes of the search tree are computed via column generation.
Lower bounds are then strengthened by the dynamic generation of valid inequalities. In Section \ref{sec:columnGeneration},
we present the proposed column generation algorithm for solving the linear relaxation of 
the F-VRPTW formulation \eqref{eq:costs}-\eqref{eq:binary}. In Section \ref{sec:accelerating}, 
we describe some acceleration strategies. Cutting planes and branching strategies are then discussed 
in Sections \ref{sec:cuttingPlanes} and \ref{sec:branching}, respectively. Finally, 
an alternative BPC algorithm which does not exploit the fragility-feasible route characterization 
is introduced in Section \ref{sec:alternative}. This algorithm will be used as a benchmark against the proposed BPC algorithm. 

\subsection{Column generation}\label{sec:columnGeneration}
At each node of the search tree, the BPC algorithm must
compute a lower bound by solving
the linear relaxation of \eqref{eq:costs}-\eqref{eq:binary}
possibly augmented by the constraints implied by the branching
decisions or valid inequalities previously added.
However, due to the extremely large number of variables in \eqref{eq:costs}-\eqref{eq:binary},
this linear problem cannot be solved directly and, consequently, 
an iterative column generation procedure needs to be invoked \citep[see, e.g.,][]{lubbecke+d2005}.
A column generation iteration first solves
a Restricted Master Problem (RMP) where only a small subset $\Omega' \subset \Omega$
of feasible routes is considered. Let $\bar{x}$ be the solution of the RMP.
The iteration then proceeds by solving a pricing problem to verify the optimality
of $\bar{x}$ for the whole linear relaxation. When the optimality check fails, a new set of routes
(identified by the pricing problem) is added to the RMP, and the process iterates. Otherwise, 
the column generation procedure stops as the current RMP solution yields a valid lower bound.

In Section \ref{sec:subForm}, we define the pricing problem, which turns out to be a variant of the ESPPRC.
This problem is solved using the labeling algorithm described in Section \ref{sec:labeling}.

\subsubsection{Pricing problem formulation}\label{sec:subForm}

To verify the optimality of the current RMP solution in the column generation procedure, 
the pricing problem searches for routes with a negative reduced cost. Let
$\pi_i$ for all $i\in \mcN_c$ be the dual variables associated with constraints~\eqref{eq:covering}.
The reduced cost of a route can be expressed as follows:
\begin{equation}\label{eq:reduced}
\bar{c}_r=c_r - \sum_{i \in \mcN_c}a_{ir}\pi_i.
\end{equation}
Thus, the pricing problem minimizes \eqref{eq:reduced} over the set of all
feasible routes $\Omega$.

Any route $r\in \Omega$ can be represented as a path in network $G$. By 
combining expressions \eqref{eq:routeCost} and \eqref{eq:reduced},
it is possible to express the reduced cost $\bar{c}_r$ of a route 
$r=(v_0, v_1,\ldots, v_p, v_{p+1})$ as the sum
of the contributions of its arcs:
$ \displaystyle \bar{c}_r=\sum_{i=1}^{p}\bar{c}_{v_i,v_{i+1}}$, where
\begin{equation}
\label{eq4}
\bar{c}_{v_i,v_{i+1}} = \left\{ 
\begin{array}{ll}
c_{v_i,v_{i+1}} & \textrm{if }v_i = 0,\\
c_{v_i,v_{i+1}} - \pi_{v_i}&\textrm{otherwise.} \\
\end{array}
\right.
\end{equation}
This property can be suitably used to express the column generation pricing problem
as an arc-flow formulation which is presented in Appendix~\ref{app:pp}. This 
formulation is, however, a non-linear integer program that is hard to solve via 
state-of-the-art solvers. On the other hand, we observe that the pricing problem 
corresponds to an ESPPRC, where the resources ensure route feasibility with respect 
to the time windows, the vehicle capacity, and the fragility constraint. 

\subsubsection{Labeling algorithm}\label{sec:labeling}

In the vehicle routing literature, the ESPPRCs are frequently solved by dynamic
programming, which is usually implemented by means of a labeling algorithm \citep{irnich+d2005}.
In our specific problem setting, given that the fragility constraint is dealt with in the
pricing problem, the labeling algorithm is potentially very challenging. In fact, given
a route visiting a set of customers in a specific order, there are many possible
ways to position the freight in the vehicle. Although some of these configurations
may be infeasible, there could still exist a large number of feasible loading configurations.
A potential labeling algorithm would create one new label for each feasible configuration or, if clever modeling is used to avoid symmetry between the identical stacks as in \citet{cherkesly+dil2016}, 
one new label for each feasible stack-anonymous configuration.  
However, both intuition and preliminary computational tests for the simplest case
with stacks of height $k=2$ suggest that this method is not efficient due to
the large number of generated labels \citep{altman2017}.

We then directed our research efforts towards the theoretical results of 
Propositions \ref{prop:characterization} and \ref{prop:simplerChar}. In particular,
the characterization of a fragility-feasible route enables us to carry 
information about the \emph{existence} of a
feasible loading configuration for a given partial route. In practice, as long as at least one feasible configuration exists, our labeling algorithm extends partial routes. However, no attempt is made to build a specific configuration, which can be trivially retrieved a posteriori.

A labeling algorithm \cite[see][]{irnich+d2005} uses labels to represent partial paths (routes) in a network.
Paths are enumerated by extending recursively an initial label $E_0$ from the source node towards the destination. Labels are extended according to extension functions. A dominance rule is applied to eliminate partial routes that cannot yield complete optimal routes. These algorithmic components are described next. 

\bigskip
\noindent {\emph{\bf Label definition}}.
A partial route $r=(0,\ldots, i)$, $i\in \mcN \cup \{ 0 \}$, is encoded by
a label $E_r=(C_r,[V_r^l]_{l\in N_c}, T_r,L_r,A^\mcH_r,A^\mcL_r,X_r)$ with the following $n+6$ components:
\begin{itemize}
\item one component $C_r$ accounting for the reduced cost;
\item $n$ binary components $V_r^l$, $l \in \mcN_c$, indicating whether or not 
customer node $l$ is \emph{unreachable}, i.e., it cannot be visited in any feasible extension of 
route $r$ because it has already been visited or because its time window or the vehicle capacity cannot be satisfied; 
\item one component $T_r$ indicating the (earliest) service starting time at node $r$;
\item one component $L_r$ accounting for total collected load;
\item two components $A^\mcH_r$ and $A^\mcL_r$ providing the total number of non-fragile and fragile
items collected, respectively;
\item one component $X_r$ indicating the maximum number of non-fragile items that
can still be collected according to Propositions \ref{prop:characterization} and \ref{prop:simplerChar}.
\end{itemize}
The initial label at node 0 is $E_{r_0}=(0,[0]_{l\in \mcN_c}, 0,0,0,0,Q)$, where partial route $r_0 = (0)$.

\bigskip
\noindent {\emph{\bf Label extension functions}}. Given a partial route $r = (0, \ldots, i)$, $i\in \mcN \setminus \{ n+1 \}$, with label $E_r=(C_r,[V_r^l]_{l\in N_c},T_r,L_r,A^\mcH_r,A^\mcL_r,X_r)$, it can be 
extended along an arc $(i,j)\in \mcA$ using the following label extension functions to 
yield a new partial route $r' = (0, \ldots, i,j)$ represented by the label 
$E_{r'}=(C_{r'},[V_{r'}^l]_{l\in \mcN_c},T_{r'},L_{r'},A^\mcH_{r'},A^\mcL_{r'},X_r)$: 
\begin{align}
C_{r'}  &=C_r+\bar c_{ij} \label{eq:costExt}\\
T_{r'}  &= \left\{ 
\begin{array}{ll} \label{eq:timeExt}
\max\{T_r+t_{ij}, \uw_j\} &  \mbox{if $j \neq n+1$} \\
T_r+t_{ij} &  \mbox{otherwise} 
\end{array} \right. \\
L_{r'}  &= L_r+q_j \label{eq:loadExt} \\
A_{r'}^\mcH&=\left\{
\begin{array}{ll}
A_r^\mcH+q_j & \mbox{if $j \in \mcH$} \\
A_r^\mcH & \mbox{otherwise}
\end{array} \right. \\
A_{r'}^\mcL&=\left\{
\begin{array}{ll}
A_r^\mcL + q_j & \mbox{if $j \in \mcL$} \\
A_r^\mcL & \mbox{otherwise} 
\end{array} \right. \\
V_{r'}^l&=\left\{
\begin{array}{ll}\label{eq:elemExt}
V_r^l +1 & \mbox{if $j=l$ } \\
V_r^l & \mbox{if $j = n+1$ } \\
\max \{V_r^l, {\cal Z}_{jl}(T_{r'},L_{r'})\} & \mbox{otherwise},
\end{array} \right. 
\qquad\qquad \forall l \in \mcN_c \\
X_{r'} & = \left\{
\begin{array}{ll}\label{eq:itapExt}
\min\{X_r, Q-L_{r'}\} & \mbox{if $j \in \mcL$ and $A^\mcH_{r'}+F(A^\mcL_{r'}) \geq k $} \\
\min\{X_r, U_{r'}\} & \mbox{if $j \in \mcL$ and $A^\mcH_{r'}+F(A^\mcL_{r'}) < k $} \\
X_{r} - q_j & \mbox{otherwise,}
\end{array} \right.
\end{align}
where $q_{n+1} = 0$, ${\cal Z}_{jl}(T_{r'},L_{r'}) = 1$ if $T_{r'} + t_{jl} > \ow_l$ or $L_{r'} + q_l > Q$ and 0 otherwise, and $U_{r'} =Q-(L_{r'} - L_{r'} \bmod k + k)$ is defined as in Proposition \ref{prop:characterization} (recall that $Q = mk$). 

Extension function \eqref{eq:itapExt} relies on Proposition~\ref{prop:characterization} and the 
vehicle capacity constraint to compute the maximum
number of non-fragile items that can still be loaded after visiting a customer.  
In the first case, given that $A^\mcH_{r'}+F(A^\mcL_{r'}) \ge k $, the loading of the $q_j$ fragile items at node $j$ does not affect the maximum number of non-fragile items
that can still be loaded in virtue of the fragility constraint. However, if the 
capacity constraint becomes binding, the term $Q - L_{r'}$ in the minimum function 
ensures that $X_{r'}$ does not exceed the residual available space (i.e., $X_{r'} \le Q-L_{r'}$). 
The second case is similar 
but, when  $A^\mcH_{r'}+F(A^\mcL_{r'}) < k$,
Proposition \ref{prop:characterization} imposes a maximum number of non-fragile items $U_{r'}$. In this case,
the capacity limit is implicit because $Q-L_{r'} > U_{r'}$.
Finally, in the third case, because $j\in \mcH$ or $j = n+1$, the extension function simply decreases 
the number of non-fragile items by $q_j$.

The obtained label $E_{r'}$ is declared feasible if $j = n+1$ or if $V_{r'}^l \leq 1$ for all $l\in \mcN_c$,
$T_{r'} \leq \ow_j$, $L_{r'}\leq Q$, and $X_{r'} \geq 0$. Otherwise, it is discarded. 
Note that all these conditions except the last one can be verified before performing the 
label extension to a node $j\in \mcN_c$ by checking the value of $V_r^j$ as follows. If $V_r^j = 1$, then at least one of these conditions will not be met and there is no need to compute $E_{r'}$. Otherwise, the 
extension must be performed and the condition $X_{r'} \geq 0$ must be checked to 
determine the feasibility of label $E_{r'}$. 

\bigskip
\noindent {\emph{\bf Dominance rule}}. Let $r_d$, $d = 1,2$, be two partial routes, both
ending at the same node $i \in \mcN$ and represented by the labels
$E_{r_d} = (C_{r_d},[V^l_{r_d}]_{l\in \mcN_c},T_{r_d},L_{r_d},A^\mcH_{r_d},A^\mcL_{r_d}, X_{r_d})$, $d = 1,2$. 
We say that $E_{r_1}$ dominates $E_{r_2}$ if
\begin{enumerate}
\item[$(c7)$] any feasible (single- or multiple-arc) extension $e$ of $r_2$ is also feasible for $r_1$, and
\item[$(c8)$] for any such extension $e$, the inequality $C_{r_1 \oplus e} \le C_{r_2\oplus e}$ holds, where 
symbol $\oplus$ denotes the concatenation operator.  \label{eq:dom2}
\end{enumerate}
The dominated label $E_{r_2}$ can then be discarded. However, when multiple labels dominate each other, one of them must be kept.

The dominance conditions $(c7)$ and $(c8)$ cannot be checked easily in practice. In general, sufficient conditions defining a so-called dominance rule are used instead. In particular, when all label extension functions are monotone (either non-decreasing or non-increasing), it is easy to provide such a dominance rule \citep[see, e.g.,][]{irnich+d2005}. 
In our case, monotonicity does not hold for extension function $\eqref{eq:itapExt}$. Consequently, the dominance rule that we introduce in the following proposition is more complex and will be proven.  
\begin{prop}\label{prop:dom} Route $r_1$ dominates route $r_2$ if the following relations hold:
\begin{align}
C_{r_1} & \leq C_{r_2} & \label{eq:costDom} \\
T_{r_1} & \leq T_{r_2} & \label{eq:timeDom} \\
L_{r_1} & \leq L_{r_2} & \label{eq:loadDom} \\
V^l_{r_1} &\leq V^l_{r_2}, & \forall l \in \mcN_c \label{eq:elemDom}\\
 \min\{X_{r_1}, Q-L_{r_1} - (k-A^\mcH_{r_1}-1)\} &\geq  Q - L_{r_2} & \label{eq:fragDom}
\end{align}
\end{prop}
In the following, we use ${\cal{R}}(\sigma_1,\sigma_2)$ as a shorthand notation for the set of conditions
\eqref{eq:costDom}-\eqref{eq:fragDom} when $r_1 = \sigma_1$ and $r_2 = \sigma_2$. 
Before proving Proposition \ref{prop:dom}, we state and prove the following preliminary lemma.

\begin{lemma}\label{prop:lemma} If relations ${\cal R}(r_1,r_2)$ are true, then relations ${\cal R}(r_1\oplus e,r_2 \oplus e)$ also hold for any extension $e$ of $r_2$.
\end{lemma}

\noindent\emph{Proof.}

We begin by showing that the statement holds true for any single-arc extension $e$ made of
a generic arc $(i,j)$. To highlight that routes $r_d$, $d = 1,2$, end at node $i$, we denote their labels by $E_{r_d} = (C_{di},[V^l_{di}]_{l\in \mcN_c},T_{di},L_{di},A^\mcH_{di},A^\mcL_{di}, X_{di})$. Similarly, the labels for routes $r_d\oplus e$ are written $E_{r_d\oplus e} = (C_{dj},[V^l_{dj}]_{l\in \mcN_c},T_{dj},L_{dj},A^\mcH_{dj},A^\mcL_{dj}, X_{dj})$.

Given that the extension functions \eqref{eq:costExt}-\eqref{eq:loadExt} and \eqref{eq:elemExt} are non-decreasing, it is trivial to show that the statement holds for relations \eqref{eq:costDom}-\eqref{eq:elemDom}. We then concentrate
on relation \eqref{eq:fragDom} for which we want to prove that:
\begin{equation}
\min\{X_{1i}, Q-L_{1i} - (k-A^\mcH_{1i}-1)\} \geq Q - L_{2i} \implies \min\{X_{1j}, Q-L_{1j} - (k-A^\mcH_{1j}-1)\} \geq Q - L_{2j}. \label{eq:implyLemma}
\end{equation}
Let us assume that $\min\{X_{1i}, Q-L_{1i} - (k-A^\mcH_{1i}-1)\} \geq Q - L_{2i}$. 
We distinguish three cases related to the extension of $r_1$, one for each case in \eqref{eq:itapExt}.
\begin{itemize}
\item If $j\in \mcH$ or $j = n+1$, then 
\begin{align*}
 \min\{X_{1j}, Q-L_{1j} - (k-A^\mcH_{1j}-1)\} 
 	&  = \min\{X_{1i} - q_j, Q-L_{1i} - (k-A^\mcH_{1i}-1)\}\\
	& \geq \min\{X_{1i}, Q-L_{1i} - (k-A^\mcH_{1i}-1)\} - q_j  \\
        &  \geq Q - L_{2i} -q_j \\
	& \geq Q-L_{2j},
\end{align*}
where the equality ensues from $X_{1j} = X_{1i} + q_j$ and $L_{1j} - A^\mcH_{1j} = L_{1i} + q_j - (A^\mcH_{1i} + q_j) = L_{1i} - A^\mcH_{1i}$, and the second inequality from the assumption. 
\item If $j\in \mcL$ and $A^\mcH_{1j}+F(A^\mcL_{1j}) \geq k$, then
\begin{align*}
\min\{X_{1j}, Q-L_{1j} - (k-A^\mcH_{1j}-1)\}  
	& = \min\{X_{1i}, Q-L_{1i} - q_j, Q-L_{1i} - (k-A^\mcH_{1i}-1) - q_j\} \\
	& \geq \min\{X_{1i}, Q-L_{1i}, Q-L_{1i} - (k-A^\mcH_{1i}-1)\} - q_j \\
	& = \min\{ \min\{X_{1i}, Q-L_{1i} - (k-A^\mcH_{1i}-1)\}, Q-L_{1i}\} - q_j \\
        & \geq \min\{Q - L_{2i}, Q-L_{1i}\} -q_j \\
	& \geq Q-L_{2j},
\end{align*}
where the equality is derived from $X_{1j} = \min \{ X_{1i}, Q-L_{1i} - q_j\}$ and $L_{1j} - A^\mcH_{1j} = L_{1i} - A^\mcH_{1i} + q_j$, the second inequality from the assumption, and the last from $L_{1i} \leq L_{2i}$.
\item If $j\in \mcL$ and $A^\mcH_{1j}+F(A^\mcL_{1j}) < k$, then 
\begin{align*}
\min\{X_{1j}, Q-L_{1j} - (k-A^\mcH_{1j}-1)\}  
	& = \min\{X_{1i}, U_j, Q-L_{1i} - (k-A^\mcH_{1i}-1) - q_j\} \\
	& \geq \min\{X_{1i}, Q-L_{1i} - (k-A^\mcH_{1i}-1) -q_j\}   \\
	& \geq \min\{X_{1i}, Q-L_{1i} - (k-A^\mcH_{1i}-1)\} - q_j  \\
	& \geq Q - L_{2i} -q_j \\
 	& \geq Q-L_{2j},
\end{align*}
where the equality arises from $X_{1j} = \min \{ X_{1i}, U_j\}$ and $L_{1j} - A^\mcH_{1j} = L_{1i} - A^\mcH_{1i} + q_j$, the first inequality from $U_j \geq Q-L_{1i} - (k-A^\mcH_{1i}-1) - q_j$ (as discussed next), and the third inequality from the assumption. Showing that $U_j= Q -L_{1i} - q_j -k+L_{1j} \bmod k \geq Q-L_{1i} - (k-A^\mcH_{1i}-1) -q_j$ is equivalent to showing that $L_{1j} \bmod k \geq A^\mcH_{1i}+1$. The latter is true because
\begin{equation*} 
L_{1j} \bmod k = (A^\mcH_{1j}+A^\mcL_{1j}) \bmod k = (A^\mcH_{1j}+F(A^\mcL_{1j})) \bmod k 
= A^\mcH_{1j}+F(A^\mcL_{1j}) \geq A^\mcH_{1i} +1,
\end{equation*}
where the second equality stems from $A^\mcL_{1j} \geq 1$ (as $j\in \mcL$), 
yielding $A^\mcL_{1j} = F(A^\mcL_{1j})$, the third equality from $A^\mcH_{1j}+F(A^\mcL_{1j}) < k$, and the inequality 
from $A^\mcH_{1j} = A^\mcH_{1i}$ and $F(A^\mcL_{1j}) \geq 1$ (as $j\in L$). 
\end{itemize}
This proves that the implication \eqref{eq:implyLemma} holds and that the statement is true for single-arc extensions. By induction on the number of single-arc extensions, the lemma is also true for any extension $e$ with an arbitrary number of arcs. $\blacksquare$

\bigskip
\noindent\emph{Proof of Proposition \ref{prop:dom}.}

To prove this proposition, we show that, if relations ${\cal R}(r_1,r_2)$ are satisfied, then conditions $(c7)$ and $(c8)$ are also satisfied. Let $e$ be a feasible extension of $r_2$ that ends at a node $j\in \mcN$. Denote by $E_{r_d\oplus e} = (C_{dj},[V^l_{dj}]_{l\in \mcN_c},T_{dj},L_{dj},A^\mcH_{dj},A^\mcL_{dj}, X_{dj})$, $d=1,2$, the labels associated with $r_d\oplus e$. Because $r_2 \oplus e$ is feasible, we have $V^l_{2j}\leq 1$ for all $l\in \mcN_c$, $T_{2j} \leq \ow_j$, $L_{2j} \leq Q$, and 
$X_{2j} \geq 0$. Assuming that relations ${\cal R}(r_1,r_2)$ hold, we deduce by Lemma \ref{prop:lemma} that 
\begin{align}
C_{1j} & \leq C_{2j} \label{eq:relC} \\
V^l_{1j} & \leq V^l_{2j}\leq 1, & \forall l\in \mcN_c \label{eq:relV} \\
T_{1j} & \leq T_{2j} \leq \ow_j \\ 
L_{1j} & \leq L_{2j} \leq Q \label{eq:relL} \\
 \min\{X_{1j}, Q-L_{1j} - (k-A^\mcH_{1j}-1)\} & \geq Q - L_{2j}. \label{eq:relX}
 \end{align}
From inequality \eqref{eq:relX}, we can easily deduce that 
\begin{equation} 
X_{1j} \geq X_{2j} \geq 0 \label{eq:relX2}
\end{equation}
because $Q - L_{2j} \geq X_{2j}$ by definition of $X_{2j}$. Therefore, 
inequalities \eqref{eq:relV}--\eqref{eq:relL} and \eqref{eq:relX2} indicate 
that $r_1\oplus e$ is feasible, i.e., that condition $(c7)$ holds, 
whereas inequality \eqref{eq:relC} implies condition $(c8)$. 
$\blacksquare$

\subsection{Acceleration strategies}\label{sec:accelerating}
Acceleration strategies play a key role for the development of efficient BPC algorithms.
In this section we describe the adopted procedures, all of them
geared towards improving the performance of the labeling algorithm in the column
generation step.

\subsubsection{Decremental state space relaxation}
An efficient technique to generate negative reduced cost elementary paths is the
decremental state space relaxation (DSSR) introduced independently by \cite{boland+dd2006}
and \cite{righini+s2008}. The pricing problem is initially solved without
the label components  $V^1,\ldots,V^n$, thus completely relaxing the elementarity
requirements. If no negative reduced cost paths are found in this way, 
some of the corresponding label components are dynamically added, and the 
process is iterated until elementary paths with a negative reduced cost are found or 
the reduced cost of the shortest path is non-negative. As proposed in \cite{desaulniers+lh2008},
at each iteration of the column generation, instead of initializing the algorithm
with an empty set of visit-label components, we use the components generated 
in the previous iteration. 
\subsubsection{The $ng$-path relaxation}\label{sec:ngPaths}
The technique described in this section consists in partially relaxing the
elementarity requirements of partial routes.
The adopted path relaxation, called $ng$-path, was introduced by \cite{baldacci+mr2011}.
With each node $i \in \mcN_c$, the approach associates a subset of nodes $N_i \subseteq \mcN_c$ 
such that $i\in N_i$ and $|N_i|\leq \Delta^0$, where $\Delta^0$ is a predefined integer parameter. 
Typically, $N_i$ contains the closest customers of $i$. 
Given a partial route $r=(v_0,\ldots,v_p)$, the subsets $N_i$ allow us to define a new subset $\Pi(r)$
whose elements are prevented to be direct extension candidates for $r$. The subset
$\Pi(r)$ is defined as $\Pi(r)=\{v_i \in r ~|~ v_i \in \bigcap_{l=i+1}^p N_{v_l},~ i = 1,\ldots,p-1\} \cup \{v_p\}$.
The value of $\Delta^0$ can influence the degree of elementarity of an $ng$-path. A
small value of $\Delta^0$ allows the $ng$-paths to contain many cycles, while
 $\Delta^0=|\mcN_c|$ imposes elementarity. Allowing cycles generally makes the corresponding
ESPPRC easier to solve, but the quality of the provided lower bound
deteriorates. Hence the right choice of $\Delta^0$ is key for the overall
computational efficiency of the BPC algorithm. 

Instead of choosing neighborhoods of the same size for all nodes, \cite{contardo+dl2015} 
propose to adjust them dynamically and individually in a DSSR fashion as follows. 
First, the size of each neighborhood is initially set to a relatively small value $\Delta^0$ 
(10 in our tests). Column generation is then applied to solve a linear relaxation. If the 
computed linear relaxation solution contains routes with a cycle, neighbors are 
added to the nodes in these cycles in such a way that they become forbidden. The 
corresponding columns are then removed from the RMP and column generation is started 
over again. The process is repeated until the linear relaxation solution is 
free of cycles or when it is not possible to forbid the remaining cycles because a 
maximum of $\Delta^+$ additional neighbors has already been added to 
a node ($\Delta^+ = 10$ for our tests). 

To implement $ng$-path, the extension function of the label components
$V^1,\ldots,V^n$ needs to be suitably modified 
\citep[see][for a detailed description of the technique]{desaulniers+mr2014,costa+cd2019}.

\subsubsection{Heuristic dynamic programming}
To rapidly obtain routes with negative reduced costs, we adopt graphs
with reduced size. In particular, at each iteration, the labeling algorithm
is first executed on a simplified graph $G'$ containing only a subset $\mcA'$ of the
arcs in $\mcA$.
In case no route with negative reduced cost is found, the labeling algorithm
is run again, but on the  complete graph. As it is detailed in
\cite{desaulniers+lh2008}, the choice of the initial set of arcs $\mcA'$
depends on the current modified costs of the arcs \eqref{eq4}. Thus, the set $\mcA'$
changes at every iteration, according to the dual values
corresponding to RMP solution. More precisely,
for every node $i \in \mcN$, all incoming arcs and all outgoing arcs are sorted according to their modified cost in
two separated lists $I_i$ and $O_i$.
If an arc $(i,j)$ is such that its position in both $I_j$ and $O_i$ is larger than a given threshold $\tau$, then
arc $(i,j)$ is eliminated. In this study, we set $\tau=10$. 

\subsection{Cutting planes}\label{sec:cuttingPlanes}
At a given node of the branch-and-bound tree, the column generation process provides
a linear relaxation solution that may be fractional, as well as a valid lower bound on 
the value of the node. When this solution is fractional, it is possible to strengthen the obtained lower bound by looking for potential violated valid inequalities.

We consider two families of inequalities. The first was
introduced in \cite{jepsen+psp2008} and is known as \emph{subset-row} inequalities,
which are  Chv{\'a}tal--Gomory rank-one cuts defined over subsets of the set partitioning constraints~\eqref{eq:covering}.
In the general case, subset-row inequalities are expressed as:
\begin{align}
\sum_{r\in\Omega} {\Big\lfloor} \frac{1}{k}\sum_{i\in W}{a_{ir}}{\Big\rfloor} x_r \leq {\Big\lfloor} \frac{|W|}{k} {\Big\rfloor}, ~~ \forall W \subseteq \mcN_c, ~ 2\leq k \leq |W|.\nonumber
\end{align}
Following \cite{jepsen+psp2008}, we only consider cuts obtained with
$|W|=3$ and $k=2$, resulting in
\begin{align}
\sum_{r\in{\Omega_W}} x_r \leq 1, ~~ \forall W\subseteq \mcN_c~: |W|=3,
\end{align}
where ${\Omega_W}$ is the subset of paths visiting at least two customers in $W$.
In a column generation method, the addition of subset-row inequalities to the RMP
requires several adjustments to the pricing problem \citep[see, e.g.,][]{desaulniers+ds2011}, 
as well as careful management when the number of added cuts increases.
The implementation proposed in this study, follows the one described in \cite{desaulniers+lh2008}. 
Note that, given the length of the routes generated for our tests, it was not worthwhile 
to apply a more sophisticated variant of the subset-row cuts such as the 
limited-node-memory subset-row cuts \citep{pecin+ppu2017}. 

The second family of inequalities we consider is the 2-path cuts,
introduced by
\citet{kohl+dmss1999} for the VRPTW. Let
$W\subseteq \mcN_c$ be a subset of nodes such that
it is not possible (e.g., because of the time windows or the vehicle capacity) 
to serve all customers in
$W$ by a single route.
The corresponding 2-path inequality is given by 
\begin{align}
\sum_{r\in \Omega} \mu^W_r x_r \ge 2, 
\end{align}
where $\mu^W_r$ is the number of times route $r$ ``enters''
the set $W$, i.e., $\mu^W_r=\sum_{i \in \mcN \setminus W}\sum_{i \in W} \eta_{ijr}$, where
non negative parameter $\eta_{ijr}$ is equal to the number of times route $r$ traverses arc $(i,j)\in \mcA$. 
We separate these cuts by a heuristic procedure similar to that of 
\cite{desaulniers+lh2008}. In particular, we enumerate all
subsets $W\subseteq \mcN_c$ such that 1) $|W|$ is less than or equal to a
given value (15 in our studies), 2)
the flow entering $W$ is strictly less than two,
and 3) the nodes in $W$ are connected in the support graph of the current
linear relaxation solution.
If $\sum_{i\in W} q_i > Q$, then a
violated inequality is found.
Otherwise,  we use dynamic programming to solve a minimum Hamiltonian path problem with time windows over the set $W \cup\{o,d\}$.
If the problem is infeasible, a violated inequality is found.
The new violated inequalities are added to the RMP. Furthermore,  the dual variable
associated with the 2-path cut for subset $W$ must be subtracted from the 
modified arc cost $\bar
c_{ij}$ for all arcs $(i,j)\in \mcA$ with $i\in \mcN \setminus W$ and $j\in W$.

\subsection{Branching strategies}\label{sec:branching}
To enforce integrality, we apply the following branching scheme whenever required.
First, we branch on
the total number of vehicles used. If this number is integer, we branch
on the arc-flow variables.
In this case, we select the arc $(i,j)$ with flow closest to $0.5$.
To fix the flow on this arc to 0, we simply remove $(i,j)$ from set $\mcA$. To fix it to $1$, we remove from $\mcA$ all arcs $(i,\ell), \ell \ne j$ if $i \ne 0$ and all arcs $(\ell, j), \ell \ne i$ if $j\ne n+1$.
The columns in the RMP containing any removed arc are deleted. 
Finally, the branch-and-bound tree is explored using a best-first strategy.

\subsection{Alternative algorithm}\label{sec:alternative}

To assess the effectiveness of the proposed BPC algorithm (denoted BPC-FC-PP 
because it handles the fragility constraint in the pricing problem), we implemented another BPC 
algorithm that does not exploit the fragility-feasible route characterization developed in 
Section~\ref{sec:char}. In fact, one may suspect that fragility-infeasible routes are relatively rare 
and should not often be encountered during the solution process. Therefore, it might not be 
worthwhile to consider the fragility constraint while solving the pricing problem and 
infeasible path cuts should rather be added to the linear relaxation whenever a linear relaxation solution 
contains a fragility-infeasible route. Infeasible path cuts have been used in the vehicle routing literature 
to handle difficult-to-model constraints pertaining to individual routes 
\citep[see, e.g.,][]{cordeau2006,cherkesly+dil2016}. Given an infeasible route $r^*$ composed of 
a subset of arcs $\mcA_{r^*}$, the corresponding infeasible path cut writes as: 
\begin{equation}
\sum_{r\in \Omega'} \sum_{(i,j)\in \mcA_{r^*}} \eta_{ijr} x_r \leq |\mcA_{r^*}| - 1, \label{eq:IPC}
\end{equation}
where $\Omega'$ is the set of routes that are feasible with respect to the time windows and 
the vehicle capacity and, for each arc $(i,j)\in \mcA$, non negative parameter $\eta_{ijr}$ 
is again equal to the number of times route $r\in \Omega'$ traverses arc $(i,j)$. 

In summary, the alternative BPC algorithm is identical to the BPC algorithm described above except that 
1) the labeling algorithm does not consider the $A^\mcH$, $A^\mcL$ and $X$ label components, nor the  
condition \eqref{eq:fragDom} in the dominance rule; and 2) infeasible path cuts of the 
form \eqref{eq:IPC} are generated as needed whenever a linear relaxation solution contains 
a fragility-infeasible route. To highlight that the fragility constraint is treated in the master problem, 
this BPC algorithm is denoted BPC-FC-MP. 

\section{Computational experiments}\label{sec:experiments}

In this section, we present the results of the computational experiments that we conducted to assess the effectiveness of the proposed BPC algorithm and to compare the solutions of the F-VRPTW obtained when varying the number of stacks and against the VRPTW solutions. In Section \ref{sec:instances}, we present the instances used for these tests. In Section \ref{sec:campaign}, we describe the tests performed and which algorithms were involved in these tests. In Section \ref{sec:results}, we report the computational results obtained and discuss computational performance. Finally, a comparison between the costs of the F-VRPTW and VRPTW solutions is realized in Section \ref{sec:cost}.

\subsection{Test instances}\label{sec:instances}

To perform our tests, we first selected the 12 instances in the class R1 of the well-known Solomon's VRPTW benchmark instances (that can be found at http://w.cba.neu.edu/$\sim$msolomon/r101.htm). In these instances involving 100 customers each, the customers are randomly located in a square grid and the time windows are relatively narrow with respect to the traveling times. These instances differ only by their customer time windows. Because the fragility constraint has more chances to be binding if vehicle capacity is relatively tight, i.e., if the filling rate of the vehicles is close to 100\%, and the routes do not contain too many customers in $\mcH$ and $\mcL$ intertwined (see the end of Section \ref{sec:char}), we have considered different tight vehicle capacities, namely, $Q = 48, 60, 72$. To study the impact of the number of stacks on the computational times and solution costs, we have also considered different stack heights relevant in practice, namely, $k = 2, 3, 4, 6$, yielding varying numbers of stacks $m = Q/k$ for each tested vehicle capacity $Q$.

Finally, for each instance, we assumed that a proportion $\rho$ of the customers belong to $\mcH$ (with non-fragile items) and, thus, $(1-\rho)$ belong to set $\mcL$ (with fragile items). We considered three different $\rho$ values, namely, 25\%, 50\%, and 75\%. To determine to which set each customer in $\mcN_c$ belongs to, we used the following procedure. Let us assume that the customers are numbered from 1 to $n$ according to the order provided by Solomon. Then, if $\rho = 50\%$, the customers with a number equal to $2i$, $i = 1, 2, 3, \ldots$ belongs to $\mcH$; if $\rho = 25\%$, the customers with a number equal to $4i - 3$, $i = 1, 2, 3, \ldots$ belongs to $\mcH$; and, if $\rho = 75\%$, the customers with a number equal to $4i - 1$, $i = 1, 2, 3, \ldots$ belongs to $\mcL$. 

Given the difficulty to solve some of these 100 customers instances to optimality, we created smaller instances by extracting the first 50 customers or the first 75 customers of each 100-customer instances. Thus, for $\rho = 50\%$, we ran tests on a total of $12\times 3 \times 4 \times 3 = 432$ instances. For $\rho = 25\%$ and $\rho = 75\%$, we also performed tests on instances with 75 customers, $Q = 60$ and $k = 2, 3, 4, 6$. 

\subsection{Computational campaign and environment}\label{sec:campaign}

Our computational experiments proceeded as follows. First, to assess the efficiency of algorithm BPC-FC-PP, we solved the 432 instances with $\rho = 50\%$ using both BPC algorithms, i.e., BPC-FC-PP and BPC-FC-MP. Then, to evaluate the impact of considering the fragility constraint on the computational times and solution costs, we also solved the corresponding VRPTW solutions using a classic BPC algorithm, i.e., algorithm BPC-FC-MP without generating any infeasible path cuts. Note that all these algorithms do not include all state-of-the-art strategies, such as route enumeration, variable fixing, or limited-memory Chv\'atal-Gomory rank-one inequalities \citep[see][for details]{costa+cd2019}. Therefore, the computational times achieved for the VRPTW instances are not competitive with those obtained by the state-of-the-art algorithms of \cite{pecin+cdu2017} and \cite{sadykov+up2020}. 

The results collected from these experiments allow to perform sensitivity analyses with respect to the number of customers, the vehicle capacity, and the stack height, but not with respect to the proportion of customers in sets $\mcL$ and $\mcH$. For the latter, we also solved the F-VRPTW instances with $\rho = 25\%$ and $\rho = 75\%$ (for 75 customers, $Q=60$ and $k =2,3,4,6$) using algorithm BPC-FC-PP only. 

All algorithms were coded in C/C{\scriptsize ++} using the Gencol 4.5 library and the primal simplex algorithm of Cplex 12.6.3 for solving the RMPs. Our tests were performed on a PC equipped with eight Intel Core i7-4770 processors clocked at 3.40 GHz, and 16 Gb of RAM. For all runs, a single core was used and a 2-hour time limit was imposed for solving each instance. 

\subsection{Computational performance comparison}\label{sec:results}

We start by presenting the results of our first experiment where, as mentioned above, 
we solved 432 instances using three algorithms. The detailed results for each individual instance 
can be found in Appendix~\ref{app:results}. We report in Table~\ref{tab:resTimes} average 
results computed over groups of instances with the same number of customers, vehicle 
capacity, and stack height. The results in this table are displayed in three blocks of columns. 
The first block presents for each F-VRPTW instance group, the results obtained by algorithm BPC-FC-PP, 
namely: the total number of instances (out of 12) solved to optimality within the time limit (\#Opt); the 
average computational time in seconds (T); the average integrality gap computed as $(z_{IP} - z_{LP})/z_{IP}$, 
where $z_{IP}$ and $z_{LP}$ are the optimal values of the problem and the root linear relaxation, respectively (Gap); 
and the average number of nodes explored in the search tree (\#Nodes). The second block provides   
the results obtained by algorithm BPC-FC-MP: the first four columns give the same information as for BPC-FC-PP; 
the last column indicates the average number of infeasible path cuts \eqref{eq:IPC} generated (\#IPCs). The third 
block presents for the VRPTW version of the instances, the total number of instances solved to optimality (\#Opt) 
and the average computational time (T). Note that all averages are computed over the instances that were 
solved to optimality within the 2-hour time limit by the corresponding algorithm. Finally, for each number of 
customers ($n = 50, 75, 100$), Table~\ref{tab:resTimes} contains a row giving totals and averages 
over all instances with the same number of customers.

\begin{table}[p]
\begin{center}
\scalebox{0.75}{
\begin{tabular}{c|c|c|rrrr|rrrrr|rr}
\toprule
&&& \multicolumn{9}{c|}{F-VRPTW} & \multicolumn{2}{c}{VRPTW} \\
\cline{4-14}
&&& \multicolumn{4}{c|}{BPC-FC-PP} & \multicolumn{5}{c|}{BPC-FC-MP} \\
 $n$ & $Q$ & $k$ &  \#Opt & T (s) & Gap (\%) & \#Nodes   &  \#Opt & T (s) &   Gap (\%) &  \#Nodes &  \#IPCs & \#Opt & T (s) \\
\midrule
\multirow{12}{*}{50} & \multirow{4}{*}{48} & 2 & 12 & 2 & 1.02 & 11 & 12 & 14 & 1.46 & 166 & 7 & 12 & 1 \\
& & 3 & 12 & 2 & 1.13 & 18 & 12 & 24 & 1.70 & 324 & 11 & 12 & 1 \\
& & 4 & 12 & 6 & 1.07 & 43 & 12 & 237 & 2.05 & 2827 & 252 & 12 & 1 \\
& & 6 & 12 & 3 & 1.08 & 16 & 12 & 109 & 1.81 & 1536 & 105 & 12 & 1 \\ 
\cmidrule{2-14}
& \multirow{4}{*}{60} & 2 & 12 & 15 & 1.42 & 65 & 12 & 13 & 1.42 & 78 & 0.2 & 12 & 13  \\
& & 3 &12 & 16 & 1.42 & 69 & 12 & 17 & 1.45 & 96 & 1 & 12 & 13  \\
& & 4 & 12 & 15 & 1.39 & 61 & 12 & 15 & 1.44 & 82 & 1 & 12 & 13  \\
& & 6 & 12 & 15 & 1.42 & 60 & 12 & 20 & 1.55 & 121 & 5 & 12 & 13  \\ 
\cmidrule{2-14}
& \multirow{4}{*}{72} & 2 & 12 & 11 & 1.43 & 25 & 12 & 10 & 1.49 & 35 & 0.9 & 12 & 9  \\
& & 3 & 12 & 10 & 1.40 & 23 & 12 & 9 & 1.48 & 29 & 0.2 & 12 & 9  \\
& & 4 & 12 & 10 & 1.42 & 24 & 12 & 11 & 1.49 & 35 & 0.9 & 12 & 9 \\
& & 6 & 12 & 10 & 1.39 & 26 & 12 & 10 & 1.49 & 35 & 1 & 12 & 9 \\ 
\cmidrule{2-14}
& \multicolumn{2}{r|}{Tot/Avg} & 144 & 9 & 1.29 & 35 & 144 & 40 & 1.56 & 447 & 32 & 144 & 8 \\
\midrule
\multirow{12}{*}{75} & \multirow{4}{*}{48} & 2 & 12 & 106 & 0.71 & 462 & 12 & 508 & 0.83 & 4402 & 218 & 12 & 153 \\
& & 3 & 12 & 302 & 0.63 & 1031 & 10 & 1046 & 0.75 & 7388 & 470 & 12 & 153  \\
& & 4 & 12 & 564 & 0.68 & 1880 & 8 & 1657 & 0.92 & 10942 & 809 & 12 & 153 \\
& & 6 & 12 & 170 & 0.72 & 574 & 4 & 1644 & 0.98 & 9556 & 959 & 12 & 153 \\ 
\cmidrule{2-14}
& \multirow{4}{*}{60} & 2 & 12 & 198 & 1.08 & 344 & 12 & 232 & 1.21 & 548 & 8 & 12 & 211  \\
& & 3 & 12 & 220 & 1.08 & 397 & 12 & 287 & 1.22 & 747 & 33 & 12 & 211  \\
& & 4 & 12 & 237 & 1.10 & 425 & 11 & 870 & 1.33 & 3177 & 140 & 12 & 211  \\
& & 6 & 12 & 369 & 1.20 & 572 & 8 & 1427 & 1.48 & 4292 & 267 & 12 & 211  \\ 
\cmidrule{2-14}
& \multirow{4}{*}{72} & 2 & 9 & 945 & 1.40 & 1126 & 8 & 808 & 1.41 & 1221 & 19 & 9 & 497  \\
& & 3 & 9 & 1279 & 1.43 & 1420 & 8 & 1403 & 1.50 & 2127 & 40 & 9 & 497  \\
& & 4 & 9 & 1159 & 1.49 & 1415 & 7 & 1504 & 1.52 & 2872 & 130 & 9 & 497 \\
& & 6 & 9 & 701 & 1.39 & 839 & 7 & 763 & 1.48 & 1433 & 55 & 9 & 497 \\ 
\cmidrule{2-14}
& \multicolumn{2}{r|}{Tot/Avg} & 132 & 475 & 1.04 & 844 & 107 & 908 & 1.20 & 3684 & 220 & 132 & 268 \\
\midrule
\multirow{12}{*}{100} & \multirow{4}{*}{48} & 2 & 11 & 668 & 0.63 & 2100 & 8 & 2440 & 0.87 & 14795 & 797 & 11 & 895  \\
& & 3 & 10 & 614 & 0.62 & 1926 & 8 & 1877 & 0.86 & 9082 & 522 & 11 & 895 \\
& & 4 & 11 & 1422 & 0.57 & 4248 & 5 & 4160 & 1.02 & 26056 & 2370 & 11 & 895 \\
& & 6 & 12 & 2152 & 0.57 & 5025 & 2 & 3387 & 1.04 & 22433 & 1716 & 11 & 895 \\ 
\cmidrule{2-14}
& \multirow{4}{*}{60} & 2 & 8 & 1187 & 0.97 & 2366 & 7 & 2264 & 1.06 & 5669 & 75 & 8 & 802 \\
& & 3 & 8 & 1712 & 1.05 & 2035 & 4 & 2603 & 1.21 & 9312 & 441 & 8 & 802 \\
& & 4 & 8 & 2961 & 1.10 & 4137 & 2 & 3516 & 1.07 & 4513 & 187 & 8 & 802 \\
& & 6 & 9 & 2260 & 1.01 & 2295 & 1 & 88 & 1.04 & 859 & 94 & 8 & 802 \\ 
\cmidrule{2-14}
& \multirow{4}{*}{72} & 2 & 9 & 1678 & 1.18 & 1639 & 5 & 2862 & 1.28 & 4146 & 75 & 10 & 1044\\
& & 3 & 10 & 1813 & 1.16 & 1386 & 4 & 1804 & 1.28 & 2792 & 51 & 10 & 1044 \\
& & 4 & 9 & 1335 & 1.14 & 1156 & 2 & 645 & 1.27 &2731 & 100 & 10 & 1044 \\
& & 6 & 10 & 1824 & 1.12 & 1080 & 2 & 566 & 1.27 & 2674 & 106 & 10 & 1044 \\ 
\cmidrule{2-14}
& \multicolumn{2}{r|}{Tot/Avg} & 115 & 1614 & 0.90 & 2506 & 50 & 2389 & 1.06 & 9956 & 592 & 116 & 230 \\
\bottomrule
\end{tabular}
}
\caption{Average results by number of customers, capacity and stack height ($\rho = 50\%$)}
\label{tab:resTimes} 
\end{center}
\end{table}

To assess the efficiency of BPC-FC-PP, let us compare the results obtained by 
BPC-FC-PP and BPC-FC-MP. First, we observe that both algorithms can solve 
all 50-customer instances within the time limit but BPC-FC-PP requires significantly 
less time on average. As the number of customers increases, 
BPC-FC-PP can solve a larger proportion of instances to optimality than 
BPC-FC-MP. In fact, the former succeeded to solve 115 instances with 100 customers, 
compared to only 50 instances for the latter. Thus, we can say that BPC-FC-PP clearly 
outperforms BPC-FC-MP. This is obviously explained by how the fragility constraint is treated 
in each algorithm because this is the only difference between these two algorithms. 
Handling it in the pricing problem increases the difficulty of solving this pricing 
problem but, as the results show, it yields substantially smaller average integrality gaps (note 
that the gap differences would probably be much larger for the instances with 75 and 100 customers 
if the same number of instances were solved to optimality by both algorithms) and, therefore, less 
branch-and-bound nodes. Also, the algorithm BPC-FC-MP spends additional time 
generating infeasible path cuts. In this regard, observe that, on average, more of these 
cuts are generated for smaller vehicle capacities 
showing, as expected, that the fragility constraint is more constraining for small-capacitated 
vehicles. The same is observed, in general, for larger stack heights although this might be less obvious 
when the number of optimally-solved instances drops as $k$ increases. Finally, when comparing the 
BPC-FC-PP results with those obtained for the VRPTW, we remark that about the same 
number of instances are solved to optimality (in total, 115 versus 116) but that the VRPTW 
instances are easier to solve (1614 versus 230 seconds). This is not surprising given that 
taking the fragility constraint into account in BPC-FC-PP increases the 
difficulty of solving the pricing problem. 

Next, we study separately the impact of the vehicle capacity and the stack height 
on the computational performance of algorithm BPC-FC-PP. Table \ref{tab:resTimesQ} 
reports the average results obtained by BPC-FC-PP by group of 144 instances with 
the same vehicle capacity (still with $\rho = 50\%$). For each instance group, the 
same statistics as in Table~\ref{tab:resTimes} are provided. These results (in particular, 
the total number of instances solved to optimality) indicate  
that the difficulty of solving the F-VRPTW increases with the vehicle capacity. Indeed, 
increasing vehicle capacity allows longer routes to be generated and makes the 
pricing problem harder to solve. Furthermore, as shown by the results, the 
average integrality gaps tend to increase with vehicle capacity, requiring more effort 
to derive an optimal integer solution. However, this is not corroborated by the 
average number of nodes explored that decreases 
when $Q$ increases. This may be explained by the fact that, compared to the case $Q=48$, 
less instances were solved to optimality and more subset-row cuts were generated 
for the cases $Q = 60$ and 72. Note that, from these results, it is not easy to see if the 
impact of handling the fragility constraint varies with the vehicle capacity. Nonetheless, 
given that the labeling algorithm used for the F-VRPTW requires an additional dominance 
condition \eqref{eq:fragDom} compared to the algorithm used for the VRPTW, 
we can speculate that handling the fragility constraint  
contributes to the increased difficulty of solving the pricing problem. 


\begin{table}[t]
\begin{center}
\scalebox{0.75}{
\begin{tabular}{c|rrrr}
\toprule
 $Q$ &  \#Opt & T (s) & Gap (\%) & \#Nodes   \\
\midrule
48 & 140 & 491 & 0.79 & 1413  \\
60 & 129 & 622 & 1.21 & 873  \\
72 & 122 & 826 & 1.32 & 772   \\
\bottomrule
\end{tabular}
}
\caption{Average BPC-FC-PP results by vehicle capacity ($\rho = 50\%$)}
\label{tab:resTimesQ} 
\end{center}
\end{table}

Using the same format as in Table \ref{tab:resTimesQ}, Table \ref{tab:resTimesK} displays the 
average BPC-FC-PP results for groups of 108 instances with 
the same stack height. Here, the variation  
between the results obtained for the different stack heights is 
much less pronounced. In fact, the number of solved instances 
for each stack height is the same except for $k=6$ for which 
three additional instances were solved. On the other hand, a 
slight increase of the average computational time is observed 
when $k$ increases. Given that the average integrality gaps 
are relatively constant with respect to stack height, this small 
computational time increase might be due to slightly 
harder-to-solve pricing problems. In fact, one can observe 
that, for a given route $r$ represented by label $E_r$, 
if the number 
of collected non-fragile items $A^\mcH_r$ is greater than 
or equal to $k-1$, then condition $(c2)$ of 
Proposition~\ref{prop:characterization} is always fulfilled 
for any subsequently visited customer in $\mcL$ 
and $X_r$ is often equal to $Q - L_r$. 
In this case, dominance condition~\eqref{eq:loadDom} 
becomes redundant with condition~\eqref{eq:fragDom} 
and the chances that $E_r$ dominates another label are 
increased. Thus, according to this hypothesis, 
a larger value of $k$ may reduce the number of 
dominated labels. 


\begin{table}[t]
\begin{center}
\scalebox{0.75}{
\begin{tabular}{c|rrrr}
\toprule
$k$ &  \#Opt & T (s) & Gap (\%) & \#Nodes   \\
\midrule
2 & 97 & 458 & 1.09 & 802  \\
3 & 97 & 578 & 1.10 & 831  \\
4 & 97 & 740 & 1.10 & 1362  \\
6 & 100 & 775 & 1.09 & 1143  \\
\bottomrule
\end{tabular}
}
\caption{Average BPC-FC-PP results by stack height ($\rho = 50\%$)}
\label{tab:resTimesK} 
\end{center}
\end{table}


We conducted a second series of experiments to see if 
the proportion of customers in set $\mcH$ influences the 
computational time. In this series, we solved all 75-customer instances with 
$Q=60$ and three different $\rho$ values, namely, 25\%, 50\%, and 75\%, 
using the BPC-FC-PP algorithm. Using the same statistics as above, 
the average results (over 48 instances per row) are reported in Table~\ref{tab:resTimesR}. 
These results indicate that the instances with $\rho = 50\%$ seem  
easier to solve than those with $\rho = 25\%$ and $\rho=75\%$: 
compared to the latter, the former exhibit average time reductions of 40\% and 27\%, 
respectively. It should be 
noticed that the high average time for the instances with $\rho = 25\%$ 
is mostly due to one outlier (instance R104 with $k =4$) which 
required more than 6700 seconds compared to less than 
1800 seconds for all others. Removing this outlier would bring 
down the average time to around 290 seconds for the instances with 
$\rho = 25\%$, that is, still a 12\% difference.  Consequently, 
a balanced mix of customers with fragile and 
non-fragile items tends to ease the solution process, but not 
substantially. A possible reason for explaining this tendency is that 
the fragility constraint is more binding when $\rho$ is close to 50\% and, thus, 
reduces more the solution space as discussed at the end 
of the next section.


\begin{table}[t]
\begin{center}
\scalebox{0.75}{
\begin{tabular}{c|rrrr}
\toprule
$\rho$ & \#Opt & T (s) & Gap (\%) & \#Nodes   \\
\midrule
25\% & 48 & 425 & 1.03 & 640 \\
50\% & 48 & 256 & 1.12 & 434 \\
75\% & 48 & 351 & 1.04 & 578 \\
\bottomrule
\end{tabular}
}
\caption{Average BPC-FC-PP results by $\rho$ value ($n = 75$, $Q = 60$)}
\label{tab:resTimesR} 
\end{center}
\end{table}

\subsection{Solution cost comparison}\label{sec:cost}

In this section, we compare the solutions obtained for the VRPTW and the F-VRPTW, namely, we 
compute the number of times that the optimal solution computed for the VRPTW 
is not feasible for the F-VRPTW and the induced average cost increase. As in Section \ref{sec:results}, 
we analyze the results with respect to the number of customers, the vehicle capacity, the stack height, 
and the proportion of customers in $\mcH$. The results are presented in 
Tables \ref{tab:resCostsN} to \ref{tab:resCostsR}. In all these tables, we report the following 
information by group of instances (the group definition depends on the table): the total number of 
instances for which the F-VRPTW was solved to optimality by either the BPC-FC-PP or the 
BPC-FC-MP algorithm (\#Opt); the total number of VRPTW solutions computed for these instances 
that are not feasible for the F-VRPTW (\#ModS) and, in parentheses, the percentage of instances 
with a modified solution, i.e., $100 \frac{\#ModS}{\#Opt}$; and the average cost increase (in percentage) 
induced by the modified solutions ($\Delta z_{IP}$). 

Table \ref{tab:resCostsN} provides the results with respect to the number of customers, i.e., each 
row summarizes the results for all instances with the same number of customers. From these 
results, we can clearly observe that the number of instances impacted by the fragility constraint 
increases with the number of customers, going from 44\% when $n = 50$ to 89\% when $n = 100$. 
This can be explained by the fact that the average number of routes in an optimal solution increases with 
the number of customers yielding more chances that at least one route in the solution is fragility-infeasible. 
The results also show that the average cost increase is relatively low (between 0.35\% to 0.61\%) 
and tends to decrease with the number of customers. Here, it is likely that more customers offer more 
options to find alternative fragility-feasible routes. 

\begin{table}[t]
\begin{center}
\scalebox{0.75}{
\begin{tabular}{c|rrr}
\toprule
$n$ & \#Opt & \#ModS & $\Delta z_{IP}$(\%) \\
\midrule
50 & 144 & 63 (44\%) & 0.61 \\
75 & 132 & 97 (73\%) & 0.40 \\
100 & 115 & 102 (89\%) & 0.35 \\
\bottomrule
\end{tabular}
}
\caption{Average cost increase by number of customers ($\rho = 50\%$)}
\label{tab:resCostsN} 
\end{center}
\end{table}


In Table \ref{tab:resCostsQ}, we present the results obtained when grouping the instances by 
vehicle capacity. From these results, we observe that the 
number of instances with a modified solution tends to decrease with vehicle capacity. 
This is to be expected because, as routes contain more customers, it becomes 
easier to find a feasible loading configuration. However, given that the proportion of 
instances with a modified solution seems to stall between $Q= 60$ (59\%) and $Q=72$ 
(58\%), we ran additional tests with $Q=84$. For this capacity, the proportion 
decreases to 38\%, confirming the observed trend. On the other hand, the 
results also indicate that the fragility constraint impacts less the solution 
cost as the vehicle capacity increases. This may be explained by the fact that, when 
capacity increases, the solutions contain less routes and, therefore, less of them 
are fragility-infeasible and need to be changed. 

\begin{table}[t]
\begin{center}
\scalebox{0.75}{
\begin{tabular}{c|rrr}
\toprule
$Q$ & \#Opt &  \#ModS & $\Delta z_{IP}$(\%) \\
\midrule
48 & 140 & 115 (82\%) & 0.57 \\
60 & 129 & 76 (59\%) & 0.38 \\
72 & 122 & 71 (58\%) & 0.26 \\
\bottomrule
\end{tabular}
}
\caption{Average cost increase by capacity ($\rho = 50\%$)}
\label{tab:resCostsQ} 
\end{center}
\end{table}

%
%
%
%

Next, we study the impact that handling the fragility constraint has on the cost of 
the solutions with respect to the stack height. The results for instances grouped by 
stack height are reported in Table \ref{tab:resCostsK}.  
Without surprise, they indicate that the number of instances with a modified solution 
and the average cost increase are influenced by stack height. Over all 
instances with the same height, we can observe an increase from 55\% to 
78\% of the \#ModS statistic when $k$ goes from 2 to 6, and an increase 
from 0.30\% to 0.55\% of $\Delta z_{IP}$. These increases are simply due to 
a reduction of the solution space that occurs when $k$ increases. 

\begin{table}[t]
\begin{center}
\scalebox{0.75}{
\begin{tabular}{c|rrr}
\toprule
$k$ & \#Opt &  \#ModS & $\Delta z_{IP}$(\%) \\
\midrule
2 & 97 & 53 (55\%) & 0.30 \\
3 & 97 & 61 (63\%) & 0.34 \\
4 & 97 & 70 (72\%) & 0.48 \\
6 & 100 & 78 (78\%) & 0.55 \\
\bottomrule
\end{tabular}
}
\caption{Average cost increase by stack height ($\rho = 50\%$)}
\label{tab:resCostsK} 
\end{center}
\end{table}


Finally, Table~\ref{tab:resCostsR} provides results for instances grouped 
by proportion of customers in~$\mcH$. They indicate that the fragility constraint 
is much more binding when there is a balance mix of customers with 
fragile and non-fragile items as 71\% of the optimal VRPTW solutions 
are fragility-infeasible when $\rho = 50\%$, whereas only 23\% and 33\% 
of the optimal solutions are infeasible for $\rho = 25\%$ and $\rho=75\%$, 
respectively. Furthermore, the cost increase is also larger for $\rho = 50\%$ 
as more routes in the VRPTW solutions become infeasible for the F-VRPTW. 
We explain this behavior as follows. As $\rho$ tends toward 0\%, the 
optimal solution routes tend to contain only customers in $\mcL$. They are, thus, 
fragility-feasible. Symmetrically, when $\rho$ tends toward 100\%, the routes 
tend to be composed of customers in $\mcH$ only and are, thus, also 
fragility-feasible. Therefore, the maximum number of routes that include 
both customer types and that are at risk of being fragility-infeasible 
occur when $\rho = 50\%$. 

\begin{table}[t]
\begin{center}
\scalebox{0.75}{
\begin{tabular}{c|rrr}
\toprule
$\rho$ & \#Opt &  \#ModS & $\Delta z_{IP}$(\%) \\
\midrule
25\% & 48 & 11 (23\%) & 0.32  \\
50\% & 48 & 34 (71\%) & 0.49  \\
75\% & 48 & 16 (33\%) & 0.22  \\
\bottomrule
\end{tabular}
}
\caption{Average cost increase by $\rho$ value ($n = 75$, $Q = 60$)}
\label{tab:resCostsR} 
\end{center}
\end{table}

\section{Conclusions}\label{sec:conclusions}

In this paper, we studied the F-VRPTW, a new variant of the VRPTW that considers a fragility 
constraint restricting the positioning of the items in the stacks of the vehicles. We first established 
necessary and sufficient conditions that characterize the feasibility of a route with respect to 
the fragility constraint. Then, we developed a BPC algorithm, denoted BPC-FC-PP, which 
exploits this characterization to handle efficiently the fragility constraint in the 
pricing problem. To evaluate this algorithm, we performed computational experiments on 
instances derived from some well-known VRPTW benchmark instances and that have been 
designed to yield a relatively tight fragility constraint. We compared 
the computational performance of BPC-FC-PP with that of another BPC 
algorithm, denoted BPC-FC-MP, that handles the fragility constraint in the master problem 
through infeasible path cuts. This comparison shows that handling the fragility constraint 
in the pricing problem is much more efficient than in the master problem: out of 432 instances, 
BPC-FC-PP solved 391 instances to optimality within a 2-hour time limit, whereas 
BPC-FC-MP could only solve 301 instances. Finally, we analyzed the impact 
that the fragility constraint has on the optimal solution cost with respect to the 
VRPTW solutions. Our results indicate that, for the tested instances, the 
computed VRPTW optimal solution is infeasible for the F-VRPTW in 67\% of 
the cases. This proportion increases with the number of customers and 
the stack height, but decreases with vehicle capacity. We also observed that 
the fragility constraint is at its maximal tightness when there is a balance mix 
of customers with fragile items and with non-fragile items. For the instances 
requiring a modified solution, an average cost increase of around 36\% 
is incurred. 

As possible future works, we can think about tackling more general problem variants. 
In particular, stacks of various heights could be considered as well as fragile items and 
non-fragile items that do not take the same space in a stack. 


\bibliographystyle{abbrvnat}
\bibliography{fragile}

\newpage
\appendix

\section{An arc-flow formulation for the pricing problem}
\label{app:pp}

In this appendix, we propose an arc-flow formulation for the 
column generation pricing problem. For each $(i,j)\in \mcA$,
we introduce a binary variable $y_{ij}$ that takes value 1 if the route 
traverses arc $(i,j)$, and 0 otherwise. For each node $i \in \mcN$, let
$T_i$ be a continuous variable indicating the service starting time. 
Furthermore, for each $(\gamma,\beta) \in M \times K$ potential position of an item in the vehicle and each customer $i\in \mcN_c$, we introduce a binary variable 
$z_{i\gamma \beta}$ which is equal to 1 if an item of customer $i$ is placed in position $(\gamma,\beta)$, and 0 otherwise.

The pricing problem can then be formulated as follows:
\begin{align} \label{eq5}
&\textrm{minimize} &&\sum_{(i,j) \in \mcA} \bar{c}_{ij} y_{ij}, \\
\label{eq6}
&\textrm{s.t.} \hspace{0.3cm} &&\sum_{j \in \mcN_c}y_{oj} = 1, \\
\label{eq7}
&&& \sum_{j \in \mcN} y_{ji} - \sum_{j \in \mcN} y_{ij} = 0, &\forall i \in \mcN_c,\\
\label{eq8}
&&& \sum_{j \in \mcN_c} y_{j,n+1} = 1, \\
\label{eq9}
&&& T_j \geq (T_i + t_{ij})y_{ij},  &\forall (i,j) \in \mcA,\\
\label{eq10}
&&& \uw_i \leq T_i \leq \ow_i,  &\forall i \in \mcN_c,\\
\label{eq11}
&&& \sum_{i \in \mcN_c} z_{i\gamma\beta} \leq 1, &\forall (\gamma,\beta)\in M\times K,\\
\label{eq12}
&&& \sum_{(\gamma,\beta) \in M \times K} z_{i\gamma\beta} = q_i \sum_{j \in \mcN}y_{ij}, &\forall i \in \mcN_c,\\
\label{eq13}
&&& \sum_{i \in \mcN_c} z_{i\gamma\beta} \geq \sum_{j \in \mcN_c} z_{j\gamma(\beta+1)}, &\forall (\gamma,\beta)\in M\times K\backslash\{k\},\\
\label{eq14}
&&& \sum_{i \in \mcN_c} T_i z_{i\gamma\beta} \leq \bar T + \sum_{j \in \mcN_c} (T_j-\bar T) z_{j\gamma(\beta+1)}, &\forall (\gamma,\beta)\in M\times K \backslash \{k\},\\
\label{eq15}
&&& \sum_{i \in \mcL} z_{i\gamma\beta} + \sum_{j \in \mcH} z_{j\gamma(\beta+1)} \leq 1, &\forall (\gamma,\beta)\in M\times K \backslash \{k\},\\
\label{eq16}
&&&y_{ij} \in \{0,1\}, &\forall (i,j) \in \mcA,\\
\label{eq17}
&&&z_{i\gamma\beta} \in \{0,1\}, &\forall i \in \mcN_c, \forall (\gamma,\beta) \in M \times K, \\
\label{eq17_1}
&&& T_i \geq 0 & \forall i \in \mcN,
\end{align}
where $\bar T$ is an upper bound on the service starting time at a customer and $T_0 = 0$.

Objective function \eqref{eq5} minimizes the total reduced cost of the route, where parameters $\bar{c}_{ij}$ are computed according to \eqref{eq4}.
Constraints \eqref{eq6}--\eqref{eq8} ensure that the $y_{ij}$ variables describe a path in $G$ between the origin depot node 0 and the destination depot node $n+1$. 
Constraints \eqref{eq9} impose lower bounds on the service starting times, while \eqref{eq10} enforce time window restrictions.
Relations \eqref{eq11}--\eqref{eq15} regulate the loading rules. In particular,
inequalities \eqref{eq11} ensure that no more than one item is put in a given position,
equalities \eqref{eq12} impose that the demand of each visited customer is fully satisfied,
inequalities \eqref{eq13} imply that items are stacked on the top of others, while
constraints \eqref{eq14}  ensure coherence between the order of the customer visits and the relative positions of the items in each stack.
Furthermore, relations \eqref{eq15} impose that no fragile item will be put under a non-fragile one.
Finally, \eqref{eq16}--\eqref{eq17_1} define the variable domains.

Note that the capacity constraint is not expressed explicitly but rather implicitly by the number of available positions in $M \times K$.

\section{Detailed computational results}
\label{app:results}

This appendix presents the detailed results obtained from our computational experiments. They are reported in Tables~\ref{tab:DetResults504850}--\ref{tab:DetResults756075}, where each table is dedicated to a number of customers ($n$), a capacity ($Q$), and a proportion of customers in $\mcH$ ($\rho$). For each instance and each stack height $k = 2, 3, 4, 6$, each table specifies the following information: 
\begin{description}
\item[$z_{IP}$:] The best solution cost found by algorithm BPC-FC-PP or BPC-FC-MP. It corresponds to the optimal value if one of these algorithms did not reach the time limit. If both algorithms reached it, an upper index indicates which algorithm (P for BPC-FC-PP and M for BPC-FC-MP) found the best cost;
\item[\#Veh:] The number of vehicles in the corresponding solution;
\item[T:] The total computational time in seconds (for algorithms BPC-FC-PP and BPC-FC-MP, as well as for solving the VRPTW). TL indicates that the 7200-second time limit has been reached;
\item[Gap:] The integrality gap in percentage (for algorithms BPC-FC-PP and BPC-FC-MP) if $z_{IP}$ is an optimal value. This gap is computed as $(z_{IP}-z_{LP})/z_{IP}$, where $z_{LP}$ is the lower bound achieved at the root node before adding cuts; 
\item[\#Nodes:] The total number of nodes explored in the search tree (for algorithms BPC-FC-PP and BPC-FC-MP);
\item[\#IPCs:] The total number of infeasible path cuts generated (for algorithm BPC-FC-MP); 
\item[$\Delta z_{IP}$:] The difference in percentage between the optimal values of the F-VRPTW and the VRPTW, whenever $z_{IP}$ is an optimal value. Note that all VRPTW optimal values have been computed even if some of them required to exceed the time limit.  
\end{description}

Note that Tables~\ref{tab:DetResults756025} and \ref{tab:DetResults756075} do not report any results for algorithm BPC-FC-MP. 


\begin{table}[p]
\begin{center}
\scalebox{0.75}{
\begin{tabular}{l|rr|rrr|rrrr|rr}
	\toprule
& \multicolumn{9}{c|}{F-VRPTW} & \multicolumn{2}{c}{VRPTW} \\
\cline{2-12}
&&& \multicolumn{3}{c|}{BPC-FC-PP} & \multicolumn{4}{c|}{BPC-FC-MP} \\
 Instance &  $z_{IP}$   &  \#Veh &  T (s) & Gap (\%) & \#Nodes   &  T (s) &   Gap (\%) &  \#Nodes &  \#IPCs & T (s) & $\Delta z_{IP}$(\%)	\\
\midrule
\multicolumn{12}{c}{\framebox[15mm]{$k=2$}}\\*[5pt]							
R101 & 1255.8 & 17 & 0.1 & 1.1 & 7 & 0.4 & 1.45 & 41 & 3 & 0.1 & 0.29 \\							
R102 & 1166.3 & 17 & 0.6 & 0.86 & 11 & 0.5 & 0.86 & 11 & 0 & 0.5 & 0 \\							
R103 & 1140.5 & 16 & 2.8 & 1.13 & 29 & 1.5 & 1.13 & 15 & 0 & 1.6 & 0 \\							
R104 & 1116.6 & 16 & 1.8 & 0.99 & 7 & 29.2 & 1.69 & 309 & 7 & 2 & 0.51 \\
\bottomrule
\end{tabular}
}
	\caption{Results for $n = 50$, $Q=48$ and $\rho = 50\%$}
	\label{tab:DetResults504850} 
\end{center}
\end{table}


\begin{table}[p]
\begin{center}
\scalebox{0.75}{
\begin{tabular}{l|rr|rrr|rrrr|rr}
	\toprule
& \multicolumn{9}{c|}{F-VRPTW} & \multicolumn{2}{c}{VRPTW} \\
\cline{2-12}
&&& \multicolumn{3}{c|}{BPC-FC-PP} & \multicolumn{4}{c|}{BPC-FC-MP} \\
 Instance &  $z_{IP}$   &  \#Veh &  T (s) & Gap (\%) & \#Nodes   &  T (s) &   Gap (\%) &  \#Nodes &  \#IPCs & T (s) & $\Delta z_{IP}$(\%)	\\
\midrule
\multicolumn{12}{c}{\framebox[15mm]{$k=2$}}\\*[5pt]							
R101 & 1137.1 & 14 & 0.1 & 0.5 & 5 & 0 & 0.5 & 3 & 0 & 0 & 0 \\							
R102 & 1047.1 & 13 & 1.6 & 1.78 & 19 & 1.3 & 1.78 & 19 & 0 & 1.3 & 0 \\							
R103 & 987.8 & 13 & 0.8 & 0.56 & 1 & 0.6 & 0.56 & 1 & 0 & 0.6 & 0 \\							
R104 & 957.5 & 13 & 3.9 & 1.08 & 5 & 3.9 & 1.08 & 19 & 0 & 4 & 0 \\
\bottomrule
\end{tabular}
}
	\caption{Results for $n = 50$, $Q=60$ and $\rho = 50\%$}
	\label{tab:DetResults506050} 
\end{center}
\end{table}


\begin{table}[p]
\begin{center}
\scalebox{0.75}{
\begin{tabular}{l|rr|rrr|rrrr|rr}
	\toprule
& \multicolumn{9}{c|}{F-VRPTW} & \multicolumn{2}{c}{VRPTW} \\
\cline{2-12}
&&& \multicolumn{3}{c|}{BPC-FC-PP} & \multicolumn{4}{c|}{BPC-FC-MP} \\
 Instance &  $z_{IP}$   &  \#Veh &  T (s) & Gap (\%) & \#Nodes   &  T (s) &   Gap (\%) &  \#Nodes &  \#IPCs & T (s) & $\Delta z_{IP}$(\%)	\\
\midrule
\multicolumn{12}{c}{\framebox[15mm]{$k=2$}}\\*[5pt]							
R101 & 1077 & 13 & 0 & 0.11 & 1 & 0 & 0.11 & 1 & 0 & 0 & 0 \\							
R102 & 959.8 & 12 & 0.2 & 0.36 & 1 & 0.5 & 0.69 & 9 & 2 & 0.1 & 0.23 \\							
R103 & 902.4 & 11 & 3.3 & 1.2 & 5 & 8.6 & 1.61 & 35 & 5 & 3 & 0.2 \\							
R104 & 856 & 11 & 27.2 & 2.03 & 47 & 19 & 2.03 & 45 & 0 & 19.6 & 0 \\
\bottomrule
\end{tabular}
}
	\caption{Results for $n = 50$, $Q=72$ and $\rho = 50\%$}
	\label{tab:DetResults507250} 
\end{center}
\end{table}


\begin{table}[p]
\begin{center}
\scalebox{0.75}{
\begin{tabular}{l|rr|rrr|rrrr|rr}
	\toprule
& \multicolumn{9}{c|}{F-VRPTW} & \multicolumn{2}{c}{VRPTW} \\
\cline{2-12}
&&& \multicolumn{3}{c|}{BPC-FC-PP} & \multicolumn{4}{c|}{BPC-FC-MP} \\
 Instance &  $z_{IP}$   &  \#Veh &  T (s) & Gap (\%) & \#Nodes   &  T (s) &   Gap (\%) &  \#Nodes &  \#IPCs & T (s) & $\Delta z_{IP}$(\%)	\\
\midrule
\multicolumn{12}{c}{\framebox[15mm]{$k=2$}}\\*[5pt]							
R101 & 1763.3 & 25 & 1.5 & 0.81 & 52 & 1.7 & 0.87 & 63 & 5 & 1 & 0.17 \\							
R102 & 1694.9 & 24 & 5.7 & 0.68 & 49 & 8.5 & 0.75 & 77 & 10 & 6.9 & 0 \\							
R103 & 1641.9 & 24 & 23 & 0.58 & 155 & 82.3 & 0.75 & 678 & 12 & 74.9 & 0 \\							
R104 & 1618 & 24 & 279.2 & 0.75 & 1031 & 339 & 0.76 & 1573 & 8 & 318.9 & 0 \\
\bottomrule
\end{tabular}
}
	\caption{Results for $n = 75$, $Q=48$ and $\rho = 50\%$}
	\label{tab:DetResults754850} 
\end{center}
\end{table}


\begin{table}[p]
\begin{center}
\scalebox{0.75}{
\begin{tabular}{l|rr|rrr|rrrr|rr}
	\toprule
& \multicolumn{9}{c|}{F-VRPTW} & \multicolumn{2}{c}{VRPTW} \\
\cline{2-12}
&&& \multicolumn{3}{c|}{BPC-FC-PP} & \multicolumn{4}{c|}{BPC-FC-MP} \\
 Instance &  $z_{IP}$   &  \#Veh &  T (s) & Gap (\%) & \#Nodes   &  T (s) &   Gap (\%) &  \#Nodes &  \#IPCs & T (s) & $\Delta z_{IP}$(\%)	\\
\midrule
\multicolumn{12}{c}{\framebox[15mm]{$k=2$}}\\*[5pt]							
R101 & 1569.5 & 20 & 4.5 & 1.14 & 83 & 7.9 & 1.79 & 175 & 1 & 0.7 & 1.1 \\							
R102 & 1481.3 & 19 & 25 & 0.96 & 115 & 50 & 1.58 & 265 & 2 & 8.2 & 0.7 \\							
R103 & 1384.3 & 19 & 218.5 & 1.35 & 441 & 252 & 1.39 & 615 & 8 & 170.2 & 0.04 \\							
R104 & 1350.2 & 19 & 733.3 & 1.04 & 1171 & 1094.2 & 1.04 & 2115 & 3 & 915 & 0 \\
\bottomrule
\end{tabular}
}
	\caption{Results for $n = 75$, $Q=60$ and $\rho = 50\%$}
	\label{tab:DetResults756050} 
\end{center}
\end{table}


\begin{table}[p]
\begin{center}
\scalebox{0.75}{
\begin{tabular}{l|rr|rrr|rrrr|rr}
	\toprule
& \multicolumn{9}{c|}{F-VRPTW} & \multicolumn{2}{c}{VRPTW} \\
\cline{2-12}
&&& \multicolumn{3}{c|}{BPC-FC-PP} & \multicolumn{4}{c|}{BPC-FC-MP} \\
 Instance &  $z_{IP}$   &  \#Veh &  T (s) & Gap (\%) & \#Nodes   &  T (s) &   Gap (\%) &  \#Nodes &  \#IPCs & T (s) & $\Delta z_{IP}$(\%)	\\
\midrule
\multicolumn{12}{c}{\framebox[15mm]{$k=2$}}\\*[5pt]							
R101 & 1486.2 & 19 & 0.2 & 0.63 & 7 & 0.3 & 0.68 & 13 & 1 & 0.3 & 0 \\							
R102 & 1367.6 & 17 & 3 & 1.15 & 9 & 2.1 & 1.22 & 7 & 0 & 2 & 0 \\							
R103 & 1245.9 & 16 & 4769.6 & 2.09 & 5945 & TL & 2.19 & 9154 & 151 & 442.2 & 0.31 \\
R104 & 1198$^{PM}$ & 16 & TL & - & 6350 & TL & - & 8909 & 162 & TL & - \\
\bottomrule
\end{tabular}
}
	\caption{Results for $n = 75$, $Q=72$ and $\rho = 50\%$}
	\label{tab:DetResults757250} 
\end{center}
\end{table}


\begin{table}[p]
\begin{center}
\scalebox{0.75}{
\begin{tabular}{l|rr|rrr|rrrr|rr}
	\toprule
& \multicolumn{9}{c|}{F-VRPTW} & \multicolumn{2}{c}{VRPTW} \\
\cline{2-12}
&&& \multicolumn{3}{c|}{BPC-FC-PP} & \multicolumn{4}{c|}{BPC-FC-MP} \\
 Instance &  $z_{IP}$   &  \#Veh &  T (s) & Gap (\%) & \#Nodes   &  T (s) &   Gap (\%) &  \#Nodes &  \#IPCs & T (s) & $\Delta z_{IP}$(\%)	\\
\midrule
\multicolumn{12}{c}{\framebox[15mm]{$k=2$}}\\*[5pt]							
R101 & 2165.6 & 33 & 2.1 & 0.49 & 29 & 11.7 & 0.87 & 245 & 19 & 4.8 & 0.24 \\							
R102 & 2105.4 & 32 & 260.7 & 0.83 & 1371 & 3876 & 1.29 & 25523 & 1381 & 102.2 & 0.49 \\	
R103 & 2029.2 & 32 & 55.4 & 0.69 & 133 & 5754.9 & 0.96 & 25957 & 1439 & 4952.8 & 0.02 \\	
R104 & 1992.5 & 32 & 809.7 & 0.58 & 1745 & 849.9 & 0.58 & 2688 & 130 & 487.3 & 0 \\
\bottomrule
\end{tabular}
}
\caption{Results for $n = 100$, $Q=48$ and $\rho = 50\%$}
\label{tab:DetResults1004850}
\end{center}
\end{table}


\begin{table}[p]
\begin{center}
\scalebox{0.75}{
\begin{tabular}{l|rr|rrr|rrrr|rr}
	\toprule
& \multicolumn{9}{c|}{F-VRPTW} & \multicolumn{2}{c}{VRPTW} \\
\cline{2-12}
&&& \multicolumn{3}{c|}{BPC-FC-PP} & \multicolumn{4}{c|}{BPC-FC-MP} \\
 Instance &  $z_{IP}$   &  \#Veh &  T (s) & Gap (\%) & \#Nodes   &  T (s) &   Gap (\%) &  \#Nodes &  \#IPCs & T (s) & $\Delta z_{IP}$(\%)	\\
\midrule
\multicolumn{12}{c}{\framebox[15mm]{$k=2$}}\\*[5pt]							
R101 & 1872.8 & 26 & 18.3 & 0.85 & 165 & 40 & 0.9 & 445 & 7 & 35.4 & 0.02 \\							
R102 & 1799.3 & 25 & 664.1 & 1.1 & 2013 & 4511.2 & 1.5 & 15943 & 111 & 2601 & 0.17 \\							
R103 & 1702.5 & 25 & 2934.6 & 1.12 & 3803 & 3361.6 & 1.16 & 5601 & 78 & 986.6 & 0.11 \\							
R104 & 2583.7$^{PM}$ & 43 & TL & - & 6667 & TL & - & 9211 & 78 & TL & - \\
\bottomrule
\end{tabular}
}
	\caption{Results for $n = 100$, $Q=60$ and $\rho = 50\%$}
	\label{tab:DetResults1006050} 
\end{center}
\end{table}


\begin{table}[p]
\begin{center}
\scalebox{0.75}{
\begin{tabular}{l|rr|rrr|rrrr|rr}
	\toprule
& \multicolumn{9}{c|}{F-VRPTW} & \multicolumn{2}{c}{VRPTW} \\
\cline{2-12}
&&& \multicolumn{3}{c|}{BPC-FC-PP} & \multicolumn{4}{c|}{BPC-FC-MP} \\
 Instance &  $z_{IP}$   &  \#Veh &  T (s) & Gap (\%) & \#Nodes   &  T (s) &   Gap (\%) &  \#Nodes &  \#IPCs & T (s) & $\Delta z_{IP}$(\%)	\\
\midrule
\multicolumn{12}{c}{\framebox[15mm]{$k=2$}}\\*[5pt]							
R101 & 1760.8 & 23 & 80.5 & 1.21 & 592 & 595.6 & 1.23 & 4457 & 105 & 334.9 & 0.06 \\							
R102 & 1632.2 & 22 & 371.1 & 1.15 & 603 & 613.5 & 1.3 & 1105 & 31 & 22.5 & 0.31 \\							
R103 & 1508.8 & 21 & 675.5 & 1 & 591 & 2603.4 & 1.32 & 3281 & 23 & 273.3 & 0.23 \\							
R104 & 1458.6 & 21 & 1447.5 & 0.98 & 593 & TL & 1.14 & 5476 & 119 & 1243.7 & 0.12 \\
\bottomrule
\end{tabular}
}
	\caption{Results for $n = 100$, $Q=72$ and $\rho = 50\%$}
	\label{tab:DetResults1007250} 
\end{center}
\end{table}


\begin{table}[p]
\begin{center}
\scalebox{0.75}{
\begin{tabular}{l|rr|rrr|rr}
	\toprule
& \multicolumn{5}{c|}{F-VRPTW} & \multicolumn{2}{c}{VRPTW} \\
\cline{2-8}
&&& \multicolumn{3}{c|}{BPC-FC-PP}  \\
 Instance &  $z_{IP}$   &  \#Veh &  T (s) & Gap (\%) & \#Nodes  & T (s) &  $\Delta z_{IP}$(\%)	\\
\midrule
\multicolumn{8}{c}{\framebox[15mm]{$k=2$}}\\*[5pt]							
R101 & 1552.5 & 20 & 0.9 & 0.72 & 23 & 0.7 &0 \\							
R102 & 1471 & 19 & 11.1 & 0.89 & 43 & 8.2 &0 \\							
R103 & 1388.1 & 19 & 828.5 & 1.58 & 1715 & 170.2 &0.32 \\							
R104 & 1350.2 & 19 & 1164.3 & 1.04 & 1805 & 915 &0 \\
\bottomrule
\end{tabular}
}
	\caption{Results for $n = 75$, $Q=60$ and $\rho = 25\%$}
	\label{tab:DetResults756025} 
\end{center}
\end{table}


\begin{table}[p]
\begin{center}
\scalebox{0.75}{
\begin{tabular}{l|rr|rrr|rr}
	\toprule
& \multicolumn{5}{c|}{F-VRPTW} & \multicolumn{2}{c}{VRPTW} \\
\cline{2-8}
&&& \multicolumn{3}{c|}{BPC-FC-PP}  \\
 Instance &  $z_{IP}$   &  \#Veh &  T (s) & Gap (\%) & \#Nodes  & T (s) &  $\Delta z_{IP}$(\%)	\\
\midrule
\multicolumn{8}{c}{\framebox[15mm]{$k=2$}}\\*[5pt]							
R101 & 1552.5 & 20 & 0.7 & 0.72 & 19 & 0.7 &0 \\							
R102 & 1471 & 19 & 9.6 & 0.88 & 33 & 8.2 &0 \\							
R103 & 1388.1 & 19 & 925.3 & 1.59 & 1771 & 170.2 &0.32 \\							
R104 & 1350.2 & 19 & 1158.9 & 1.01 & 1783 & 915 &0 \\
\bottomrule
\end{tabular}
}
	\caption{Results for $n = 75$, $Q=60$ and $\rho = 75\%$}
	\label{tab:DetResults756075} 
\end{center}
\end{table}

\end{document}